\documentclass[onefignum,onetabnum]{siamart171218}
\pdfoutput=1

\usepackage{amsfonts}
\usepackage{graphicx}
\usepackage{xcolor}
\usepackage[ruled, linesnumbered,algo2e]{algorithm2e}
  \DeclareGraphicsExtensions{.eps,.pdf,.png,.jpg}

\usepackage{epstopdf}

\newsiamremark{remark}{Remark}
\newsiamremark{hypothesis}{Hypothesis}
\crefname{hypothesis}{Hypothesis}{Hypotheses}
\newsiamthm{claim}{Claim}

\headers{Subspace Recycling-based Regularization Methods}{R. Ramlau and K. M. Soodhalter and Victoria Hutterer}

\title{Subspace Recycling-based Regularization Methods\thanks{DATE.
}}

\author{Ronny Ramlau\thanks{Industrial Mathematics Institute, Kepler University Linz, and Johann Radon Institute for Computational and Applied Mathematics (RICAM), Linz, Austria.({\tt ronny.ramlau@jku.at}).}
\and Kirk M. Soodhalter\thanks{	School of Mathematics, Trinity College Dublin, The University of Dublin, College Green, Dublin 2, Ireland.({\tt ksoodha@maths.tcd.ie}).}
\and Victoria Hutterer
\thanks{Industrial Mathematics Institute, Kepler University Linz, Austria. ({\tt victoria.hutterer@indmath.uni-linz.ac.at}).}
}

\usepackage{amsopn}

\def\cal{\mathcal}

\def\nm#1{\left\|#1\right\|}

\def\R{\mathbb{R}}

\def\Rn{\R^n}

\def\BA{{\bf A}}  
  
  \def\CC{{\cal C}}

  \def\CF{{\cal F}}

  \def\CK{{\cal K}}
  \def\CL{{\cal L}}

  \def\CO{{\cal O}}
  \def\CP{{\cal P}}

  \def\CS{{\cal S}}
  
  \def\CU{{\cal U}}
  \def\CV{{\cal V}}
  
  \def\CX{{\cal X}}
  \def\CY{{\cal Y}}

\def\BAs{\BA{\kern-1.5pt}}

\def\CPs{\CP{\kern-0.8pt}}

\mathcode`\@="8000
{\catcode`\@=\active \gdef@{\mkern1mu}}

\def\mydate{\number\day\ {\ifcase\month \or January\or February\or
              March\or April\or May\or June\or July\or August\or
              September\or October\or November\or December\fi}
\number\year}
\def\vek#1{\mathbf{#1}}

\def\ip#1{\left\langle#1\right\rangle}

\providecommand{\argmin}[1]{\underset{#1}{\text{{\rm argmin}}}}

\providecommand{\prn}[1]{\left(#1\right)}
\providecommand{\bigprn}[1]{\big(#1\big)}

\providecommand{\brac}[1]{\left[#1\right]}
\def\curl#1{\left\{#1\right\}}
\providecommand{\ab}[1]{\left|#1\right|}

\newcommand\restr[2]{{
  \left.\kern-\nulldelimiterspace 
  #1 
  \vphantom{\big|} 
  \right|_{#2} 
  }}

\def\Span{{\rm span}}

\def\be{\begin{equation}}
\def\ee{\end{equation}}
\def\bea{\begin{eqnarray}}
\def\eea{\end{eqnarray}}
\def\nn{\nonumber}
\def\mand{\mbox{\ \ \ and\ \ \ }}

\def\mfor{\mbox{\ \ \ for\ \ \ }}

\def\mwith{\mbox{\ \ \ with\ \ \ }}
\def\mwhere{\mbox{\ \ \ where\ \ \ }}

\def\msuchthat{\mbox{\ \ \ such that\ \ }}

\def\mselect{\mbox{\ \ \ select\ \ }}

\def\mif{\mbox{\ \ \ if\ \ \ }}

\def\motherwise{\mbox{\ \ \ otherwise}}

\def\bbmat{\begin{bmatrix}}
\def\ebmat{\end{bmatrix}}

\def\balg{\begin{algorithm}}
\def\ealg{\end{algorithm}}
\def\balgte{\begin{algorithm2e}}
\def\ealgte{\end{algorithm2e}}
\def\bthm{\begin{theorem}}
\def\ethm{\end{theorem}}
\def\blem{\begin{lemma}}
\def\elem{\end{lemma}}
\def\bprop{\begin{proposition}}
\def\eprop{\end{proposition}}
\def\bcor{\begin{corollary}}
\def\ecor{\end{corollary}}
\def\bdefin{\begin{definition}}
\def\edefin{\end{definition}}
\def\bc{\begin{cases}}
\def\ec{\end{cases}}
\newcommand\bproof[1]{\par\addvspace{1ex} \indent\textit{Proof.}\ \ #1}
\def\eproof{\hfill\cvd\linebreak\indent}
\newtheorem{exple}{Example}
\def\bex{\begin{exple}}
\def\eex{\end{exple}}
\newtheorem{exercise}{Exercise}
\def\bexer{\begin{exercise}}
\def\eexer{\end{exercise}}
\newtheorem{conjecture}{Conjecture}
\def\bconj{\begin{conjecture}}
\def\econj{\end{conjecture}}
\newtheorem{assumption}{Assumption}
\def\bass{\begin{assumption}}
\def\eass{\end{assumption}}
\newtheorem{notation}{Notation}
\def\bnot{\begin{notation}}
\def\enot{\end{notation}}
\def\brem{\begin{remark}}
\def\erem{\end{remark}}

\def\blankfootnote#1{\let\thefootnote\relax\footnotetext{#1}}
\def\cvd{~\vbox{\hrule\hbox{%
  \vrule height1.3ex\hskip0.8ex\vrule}\hrule } }

\def\bitem{\begin{item}}
\def\eitem{\end{itemize}}
\def\benum{\begin{enumerate}}
\def\eenum{\end{enumerate}}

\DeclareMathSymbol{\dprod}{\mathbin}{operators}{"3A}

\providecommand{\file}[1]{\texttt{\nolinkurl{#1}}}

\def\xpdag{x_p^\dagger}
\def\xpd{x_p^\delta}
\def\xad{x_\alpha^\delta}

\begin{document}

\maketitle
\begin{abstract}
Subspace recycling techniques have been used quite successfully for the acceleration of iterative methods for solving large-scale linear systems.  These methods often work by augmenting a solution subspace generated iteratively by a known
algorithm with a fixed subspace of vectors which are ``useful'' for solving the problem.  Often, this has the effect of inducing a
projected version of the original linear system to which the known iterative method is then applied, and this projection can
act as a deflation preconditioner, accelerating convergence.
Most often, these methods have been applied for the solution of well-posed problems.  However,
they have also begun to be considered for the solution of ill-posed problems.

In this paper, we consider subspace augmentation-type iterative schemes applied to linear ill-posed problems in a continuous
Hilbert space setting, based on a recently developed framework describing these methods.
We show that under suitable assumptions, a recycling method satisfies the formal definition
of a regularization, as long as the underlying scheme is itself a regularization.  We then develop an augmented subspace
version of the gradient descent method and demonstrate its effectiveness, both on an academic Gaussian blur model and on
problems arising from the adaptive optics community for the resolution of large sky images by ground-based extremely large
telescopes.
\end{abstract}



\section{Introduction}
{\color{black}  This paper concerns the use of subspace recycling iterative methods in the context of solving linear 
ill-posed problems.  There have already been a number of augmented or subspace recycling methods proposed
specifically to treat ill-posed problems (e.g., \cite{BR.2007,BR-2.2007,DGH.2014,dSC.2019}), 
but these methods have yet to be formally analyzed as
regularization schemes.  
More recently, an augmented scheme to accelerate the regularization of a problem in adaptive optics,
producing a large improvement in performance, was proposed \cite{ramlau2020augmented}.
However, the augmentation 
was done after the regularized problem was set up rather than actually combining these techniques.
It would thus be helpful to place all these methods into a common framework to streamline the development of new methods
and to enable them to be analyzed using a common language and framework.  In particular, it enables us in this paper to
prove under what conditions augmented approaches described by the framework are regularizations.

We take advantage of a newly-proposed general framework 
\cite{soodhalter2020survey} which can be used to describe the vast majority of these methods. This framework
allows us to perform a general analysis of augmentation and recycling schemes as regularizations.  From this, we develop
 sufficient conditions for a recycled iterative method to be a regularization, namely that
a subspace recycling scheme is a regularization if it is built from an underlying regularization method.  As a proof-of-concept,
we then leverage the framework and our analysis to propose augmented 
steepest descent and Landweber methods, which we then apply to some test problems.}

The linear inverse problem we consider is of the form
\begin{equation}\label{illposedeq}
Tx=y
\end{equation}
{\color{black} whereby we approximately reconstruct $x$ from noisy data $y^\delta$ where $y^{\delta} = y + n$, and $n$ represents the perturbation of our data which satisfies $\|n\|\le \delta$.}
 Here, $T:\CX\to\CY$ is a continuous mapping between separable Hilbert spaces $\CX , \CY$. In many relevant applications, the operator $T$ is ill-posed, i.e., \eqref{illposedeq} violates Hadamard's conditions: a solution either might not exist, it might be non-unique or it may not depend continuously on the data. In this case regularization techniques have to be employed for the stable solution of \eqref{illposedeq}, see \cite{RamlauEngl2015,Engl_Hanke_Neubauer_1996} and Section \ref{reg_techniques} for a summary on these methods.  
 
 The quality of the reconstruction for a specific problem depends heavily on available additional information on the solution, which is often directly included into the chosen reconstruction approach. E.g., source conditions of the solution are used to estimate the reconstruction error \cite{Engl_Hanke_Neubauer_1996,Louis_1989}.   Alternatively, information on the smoothness of the solution can be incorporated into Tikhonov regularization by choosing a penalty term that enforces the reconstruction to be smooth. For iterative regularization methods, the iteration can be started with a rough guess for the solution, leading generally to shorter reconstruction times and better results. Those approaches are of particular interest for inverse problems that have to be solved repeatedly in time and where the previous reconstruction can be considered as a good guess for the solution of the current problem.
This would be the case, e.g., for {\it Adaptive Optics} applications in Astronomy, where the image quality loss due to turbulences in the atmosphere is corrected by a suitable shaped deformable mirror. The shape of the deformable mirrors has to be adjusted within milliseconds. The optimal shape of the mirror is derived by solving the atmospheric tomography problem repeatedly (about every 1-2ms). Although the atmosphere changes rapidly it can be regarded {\it frozen} on a sub-millisecond timescale, meaning that the previous reconstruction of the atmospheric tomography problem can be regarded as a good approximation to the new solution \cite{YuHeRa13b,RaRo12}.

In this paper we want to consider the case that we know a finite-dimensional subspace $\CU\subset\CX$ that  
contains a significant part of the solution $x^\dagger$ which we seek. 
As we will see, this part can be easily and stably be reconstructed. What then remains is to reconstruct the part that belongs to the orthogonal complement $\CU^\perp$ of $\CU$,
possibly 
induces an acceleration in the convergence
rate of the resulting iterative method, cf. \Cref{section.basic-aug-proj}.

For finite dimensional problems, this approach is known as a {\it subspace recycling method} and is currently a hot topic in numerical linear algebra. The goal of this paper is therefore to extend the theory to ill-posed problems in an infinite dimensional setting. Specifically, we will show that the concept of subspace recycling can be used in connection with most of the standard regularization methods.

The paper is organized as follows: In Section 2 we give a short introduction to Inverse Problems and regularization methods, whereas Section 3 discusses the state of research of subspace recycling methods. In Sections 4 and 5 we show that standard linear regularization method can be combined with subspace recycling and prove that the resulting method is still a regularization. Finally, Section 6 considers recycled gradient descent methods, and Section 7 presents some numerical results.

\section{Regularization of linear ill-posed problems}\label{reg_techniques}
{\color{black}We begin with a brief review of some basic regularization theory, which we mainly base on} \cite{RamlauEngl2015,Engl_Hanke_Neubauer_1996,Louis_1989}. As mentioned above, a solution to the problem $Tx=y$ might not exist; if it exists it might not be unique and/or might not depend continuously on the data. The first two problems can be circumvented by using the generalized inverse $T^\dagger$, whereas stability is enforced by using {\it regularization methods}. To be more specific, an operator
\begin{equation}
R_\alpha :Y\to X
\end{equation}
is a regularization for a linear operator $T:X\to Y$ iff for any right hand side $y^\delta$ with $\|y-y^\delta\|\le\delta$ {\color{black} there} exists a parameter choice rule
\begin{equation}
\alpha: \R\times Y\to \R^+
\end{equation}
such that
\begin{equation}
  \lim_{\delta\to 0} R_{\alpha (\delta)}y^\delta= T^\dagger y = :x^\dagger
\end{equation}
holds.
Here, the parameter is chosen such that
$\lim_{\delta\to 0}\alpha (\delta )=0$.
Over the last decades, many regularization methods have been investigated, most prominently {\it Tikhonov regularization} and iterative methods {\color{black}like {\it Landweber iteration}~\cite{Engl_Hanke_Neubauer_1996,Landweber51}} and the {\it conjugate gradient} method. For Tikhonov regularization, the approximation $\xad$ to the generalized solution $x^\dagger$ is computed as the minimizer of the Tikhonov functional
\begin{equation}
  J_\alpha (x) = \|y^\delta -Tx\|_Y^2+\alpha \Omega (x),
\end{equation}
where $\Omega (x)$ is a suitable penalty term, e.g., $\Omega (x)=\|x\|^2$. For more general penalties we refer, e.g., to \cite{RamlauEngl2015,ScherzerBook09}. {\color{black} For Tikhonov to be rendered a regularization method, one appropriate method of choosing the regularization parameter is the discrepancy principle, where $\alpha $ is chosen such that }
\begin{equation}
  \|y^\delta - Tx_\alpha^\delta\|=\tau\delta , \hspace{1cm}\tau>1.
\end{equation}
Iterative methods are popular as their computational realization is often straightforward. Landweber iteration computes as
\begin{equation}\label{Landweber}
  x_{k+1}^\delta = x_{k}^\delta +\beta T^\ast (y^\delta -Tx_k^\delta),\hspace{1cm} 0<\beta<2/\|T\|^2.
\end{equation}
{\color{black}For iterative methods the iteration index plays the role of the regularization parameter. According to the discrepancy principle, the iteration is terminated at step $k$ when for the first time $\|y^\delta - Tx_k^\delta\|\le \tau \delta$ with $\tau > 1$ fixed. The Landweber iteration with the discrepancy principle is a regularization method~\cite{Engl_Hanke_Neubauer_1996}.} It is well known that the convergence of regularization methods for $\delta\to 0$ can be arbitrarily slow. Convergence rates can be obtained by using {\it source conditions}. Popular are H\"older type source conditions
\begin{equation}\label{sourcecond}
  x^\dagger = (T^\ast T)^\nu w
\end{equation}
{\color{black} for some $\mu \ge 0$ and $w \in X$.} A regularization method is called order optimal iff the reconstruction error can be estimated as
\begin{equation}
  \|x^\dagger - R_{\alpha (\delta)}y^\delta\| = \mathcal{O} (\delta^{\nu /(\nu+1)})
\end{equation}
provided a source condition \eqref{sourcecond} is fulfilled. {\color{black}Both Tikhonov and Landweber equipped with suitable source conditions are order optimal~\cite{Engl_Hanke_Neubauer_1996}.}

In this paper we will show that subspace recycling-based methods also form a regularization.

\section{Augmented iterative methods}\label{section.aug-iter}
Augmented iterative methods, as we discuss them here, have been proposed in the literature beginning in the mid- to late-nineties; see \cite{soodhalter2020survey} for a survey of this history.
They can be considered as a specialized version of the class of iterative methods called \emph{residual constraint} or \emph{residual projection}
methods.  We discuss these methods briefly to put them into context, but for more details, 
{\color{black} the reader should see the survey paper \cite{soodhalter2020survey}, and references therein}.
\subsection{Iterations determined by residual constraint}\label{section.constr-iter}
Many iterative methods which are commonly used to treat both well- and ill-posed linear problems
can be formulated according to a correction/constraint setup.  Our discussion here in \Cref{section.aug-iter} focuses on the
well-posed problem case for simplicity.  Let $\prn{\cdot,\cdot}_{\CX}$ and $\prn{\cdot,\cdot}_{\CY}$,
respectively, be the inner products for $\CX$ and $\CY$ and
$\prn{Tx,y}_{\CX} = \prn{x,T^{\ast}y}_{\CY}$.

Consider the case wherein we are solving
\cref{illposedeq} with no initial approximation; i.e., $x_{0}=0$.  The strategies we outline here take the abstract form of
approximating the solution $x^{\dagger}$ by $\hat{x}\in\CS$ such that the residual satisfies $y - T\hat{x} \perp \widetilde{\CS}$, where
$\CS\subset\CX$ and $\widetilde{\CS}\subset \CY$ are such that $\dim \CS = \dim \widetilde{\CS}<\infty$.
They are called the \emph{correction} and \emph{constraint} spaces, respectively.  A common example of a case for which
$\dim\CS=1$ is steepest descent (i.e., gradient descent with step size chosen to minimize the $T^{\ast}T$-norm of the error); see, e.g., \cite[Section 5.3.1]{Saad.Iter.Meth.Sparse.2003}.
\bprop\label{eqn.steepest-descent-characterization}
The $(j+1)$st step of steepest descent can be formulated as
\be\label{eqn.steepest-descent-proj}
	\mselect t_{j+1}\in \Span\curl{T^{\ast}r_{j}}\msuchthat y-T\prn{x_{j} + t_{j+1}}\perp \Span\curl{TT^{\ast}r_{j}},
\ee
where $r_{j} = y-Tx_{j}$; i.e., $\CS = \Span\curl{T^{\ast}r_{j}}$, and $\widetilde{\CS} = \Span\curl{TT^{\ast}r_{j}}$.\footnote{One could also,
with some work, use this result to
obtain a residual constraint formulation for a more general Landweber method, but it is not clear that this brings an analytic advantage.}
\eprop
\begin{proof}
Suppose we have constructed the approximation $x_{j}$ at the previous iteration.  Recall that gradient descent methods compute
$x_{j+1} = x_{j} + \alpha T^{\ast} r_{j}$, where $r_{j} = y - Tx_{j}$.  The choice of $\alpha$ (called the step size) determines the method's behavior.
Steepest descent computes
\be\nn
	\alpha_{j+1} = \argmin{\alpha\in\R}\nm{y - T(x_{j} + \alpha T^{\ast}r_{j})}_{\CY}.
\ee
Differentiating the right-hand side, setting it equal to zero, and solving for $\alpha_{j+1}$ yields the minimizer
$\alpha_{j+1} =\dfrac{\prn{r_{j},TT^{\ast}r_{j}}_{\CY}}{\prn{TT^{\ast}r_{j},TT^{\ast}r_{j}}_{\CY}}$.  One must simply show
that \cref{eqn.steepest-descent-proj} produces the same choice of $\alpha_{j+1}$.  However, this constraint condition
is equivalent to solving the equation
\be\nn
	\prn{r_{j} - \alpha_{j}TT^{\ast}r_{j}, TT^{\ast}r_{j}}_{\CY} = 0,
\ee
which when solved for $\alpha_{j}$ produces the same $T^{\ast}T$-norm-minimizing choice of $\alpha_{j}$.
\end{proof}

Steepest descent is an example of a method for which $\dim\CS=1$.  There are many methods for which the size of the
correction and constraint spaces grows at each iteration.  It is a well-known result that one can characterize the method of conjugate gradients in this way, which
we present without proof; see \cite{S1981,Saad.Iter.Meth.Sparse.2003} for proofs.
\bprop
The $(j+1)$st iterate produced by the method of conjugate gradients applied to the normal equations associated to \cref{illposedeq}
can be formulated as
\be\label{eqn.CGNE-proj}
	\mselect x_{j+1}\in \CK_{j+1}\prn{T^{\ast}T,T^{\ast}y}\msuchthat y-Tx_{j+1}\perp T\CK_{j+1}\prn{T^{\ast}T,T^{\ast}y},
\ee
where $\CK_{j+1}\prn{T^{\ast}T,T^{\ast}y} = \Span\curl{y,T^{\ast}Ty,\prn{T^{\ast}T}^{2}y,\ldots,\prn{T^{\ast}T}^{j}y}$.
\eprop

\subsection{Augmentation via constraint over a pair of spaces}
Many flavors of augmented iterative methods have been proposed to treat both well- and ill-posed problems. These methods
treat the situation in which we wish to approximate the solution using an iterative method, but we also have a fixed subspace $\CU$ which we also
want to use to build the approximation.  Using the correction/constraint formulation in \Cref{section.constr-iter}, we can describe many of these methods as
selecting approximations from a correction space $\CS_{j} = \CU + \CV_{j}\subset\CX$ with a residual constraint space
$\widetilde{\CS}_{j} = \widetilde{\CU} + \widetilde{\CV}_{j}\subset\CY$, where $\CU$ and $\widetilde{\CU}$
are spaces of fixed dimension $k$ and $\CV_{j}$ and $\widetilde{\CV}_{j}$
are generated iteratively, with dimension $j$ at step $j$.  It has been shown that this correction/constraint over the sum of subspaces reduces mathematically to applying the underlying
iterative method to a projected problem, obtaining part of the approximation from $\CV_{j}$ and then obtaining the corrections from $\CU$ afterward.  This will be explained in more detail in \Cref{section.basic-aug-proj}, but we give a brief outline here.

A projector is uniquely defined by its range and its null space.  A standard result from linear algebra shows that this is equivalent to defining the ranges of the
projector and its complement.  In this setting,
consider the projectors $P\in\CL\prn{\CX}$ and $Q\in\CL\prn{\CY}$ having the relationship that
\be\label{eqn.projector-relationship}
	TP = QT.
\ee
Let ${\rm Range}\prn{P}=\CU$ and ${\rm Null}\prn{P}=\prn{T^{\ast}\,\widetilde{\CU}}^{\perp}$, and  ${\rm Range}\prn{Q}=T\,\CU$ and ${\rm Null}\prn{Q}=\prn{\widetilde{\CU}}^{\perp}$.
One can easily show that defining the null spaces is equivalent to defining ${\rm Range}\prn{I_{\CX} - P}=\prn{T^{\ast}\,\widetilde{\CU}}^{\perp}$ and ${\rm Range}\prn{I_{\CY} - Q}=\widetilde{\CU}^{\perp}$.
For example, let $v\in\CX$ be chosen arbitrarily.  It can be decomposed as $v = v_{\mathcal{R}} + v_{\mathcal{N}}$ with $v_{\mathcal{R}}\in{\rm Range}\prn{P}$ and $v_{\mathcal{N}}\in{\rm Null}\prn{P}$.
As expected, $Pv = Pv_{\mathcal{R}} + Pv_{\mathcal{N}}=v_{\mathcal{R}}$.  We can then observe that $\prn{I_{\CX} - P}v = v_{\mathcal{R}} + v_{\mathcal{N}} - Pv_{\mathcal{R}} - Pv_{\mathcal{N}} = v_{\mathcal{N}}$.
As $v\in\CX$ was chosen arbitrarily, we see that ${\rm Range}\prn{I_{\CX} - P} = {\rm Null}\prn{P}$.  The same can be shown for $Q$.

Approximating the solution of \cref{illposedeq} by $x_{j} \in \CS_{j}$ subject to the constraint $r_{j}\perp\widetilde{\CS}_{j}$ has been shown to be equivalent to
approximating the solution $t$ to the projected problem
\be\label{eqn.proj-prob}
	\prn{I - Q}Tt = \prn{I - Q}y
\ee
by $t_{j}\in\CV_{j}$ with residual constraint space $\widetilde{\CV}_{j}$ and setting the approximation $x_{j}\in\CU+\CV_{j}$ to be 
\begin{equation}\label{eqn.x-approx-decomp}
	x_{j} = Px^{\dagger} + \prn{I-P}t_{j};
\end{equation} 
see, \cite{dSKS.2018} and \cite{Gaul.2014-phd,GGL.2013,Gutknecht.AugBiCG.2014}.

For well-posed problems, the space $\CU$ often contains vectors from previous iteration cycles applied to the current or previously-solved related problem; see, e.g., \cite{ASGC-rBiCG-Model-Red.2012,Kilmer.deSturler.tomography.2006,Parks.deSturler.GCRODR.2005,WSG.2007}.
These vectors are
often (but not always) approximate eigenvectors (of a square operator) or left singular vectors, {\color{black}guided by Krylov subspace
method convergence theory}.    For ill-posed problems, the strategy is generally to
include some known features of the image/signal being reconstructed (see, e.g., \cite{BR-2.2007,BR.2007,DGH.2014}),
 but some current work by \cite{dSC.2019} still seeks to
augment a Lanczos bidiagonalization-based solver for ill-posed, tall, skinny least-squares problems using approximate left singular
vectors to accelerate semi-convergence and reduce the influence of smaller singular values to achieve a regularizing effect.
{\color{black} See \Cref{sec.what-to-recycle} for further discussion.}

There is an alternative approach
(originally proposed in the context of domain decomposition in, e.g., \cite{Erlangga2008}) which begins with
the projectors $P$ and $Q$ (rather than treating them as a consequence of the correction/constraint formulation)
and produces the same algorithms but provides more flexibility to use methods like Landweber which do not have a clean residual constraint condition.
This will be
derived in detail in \Cref{section.basic-aug-proj}, but it is built on the simple idea that the projector $P$ can be used to get the decomposition
$x^{\dagger} = Px^{\dagger} + \prn{I-P}x^{\dagger}$.  We will show that $Px^{\dagger}$ can be cheaply computed and we can use an iterative method to
approximate $\prn{I-P}x^{\dagger}$.

\section{Linear Projection Methods}
We consider here the case that $\widetilde{\CU} = T\,\CU$.
For a finite dimensional subspace $\CU\subset\CX$  with $\dim\CU = k > 0$, we
denote the sibling subspace $\CC = T\,\CU\subset\CY$ and assume it also has dimension $k$.
We can represent these subspaces in terms of bases, which we denote with a matrix-like notation.  Let
\be\label{def_subspaces}
	\vek U = \bbmat u_{1} & u_{2} & \cdots & u_{k} \ebmat\in\CX^{k} \mand \vek C=\bbmat c_{1} & c_{2} & \cdots & c_{k} \ebmat\in\CY^{k}
\ee
be $k$-tuples of elements from $\CU$ and $\CC$, respectively, which are bases for these spaces.

\brem
To get a sense of these objects, observe that if $\CX=\CY=\Rn$ then we would have $\vek U,\vek C\in\R^{n\times k}$; so this is just the generalization of the
notion of a tall, skinny matrix.
\erem

We can consider $\vek U$ and $\vek C$ to be linear
mappings induced by a basis expansion operation, with $\vek U\in\CL\prn{\R^{k},\CX}$ and $\vek C\in\CL\prn{\R^{k},\CY}$, defined by
\be\label{UC-Def}
	\vek U\vek z = \sum_{i=1}^{k}z_{i} u_{i}\mand \vek C\vek z = \sum_{i=1}^{k}z_{i} c_{i}\mfor \vek z = \ \bbmat z_{1} & z_{2} & \cdots & z_{k} \ebmat^{T}\in\R^{k}.
\ee
Let
\be\label{vekmat}
\prn{x,\vek U}_{\CX} = \bbmat \prn{x,u_{1}}_{\CX}\\\vdots \\ \prn{x,u_{k}}_{\CX} \ebmat\in\R^{k},\mand \prn{y,\vek C}_{\CY} = \bbmat \prn{y,c_{1}}_{\CY}\\\vdots \\ \prn{y,c_{k}}_{\CY} \ebmat\in\R^{k}.
\ee
Furthermore, for $\vek M,\vek L \in\CX^{k}$, we define the bilinear mapping
\be\label{matmat}
	\prn{\cdot,\cdot}_{\CX} : \CX^{k} \times \CX^{k}\rightarrow \R^{k\times k}\mwith \prn{\vek M,\vek L}\mapsto \bigprn{\prn{m_{i},\ell_{j}}_{\CX}}_{i,j=1}^{k}\in\R^{k\times k}.
\ee
We similarly define the corresponding bilinear mapping for two elements from $\vek F,\vek G\in\CY^{k}$
\be\nn
	\prn{\cdot,\cdot}_{\CY} : \CY^{k} \times \CY^{k}\rightarrow \R^{k\times k}\mwith \prn{\vek F,\vek G}\mapsto \bigprn{\prn{f_{i},g_{j}}_{\CY}}_{i,j=1}^{k}\in\R^{k\times k}.
\ee

In particular, we have the following relation:
\blem\label{Unorm}
For the operator $\vek U$ (and similarly $\vek C$) holds
\be
\|\vek U z\|_{\CX}^2\le \|(\vek U , \vek U)\|_F \|z\|^2_{\R^k}
\ee
\elem
\begin{proof}
  With $z=(z_1,\dots , z_k)\in\R^k$, $\vek Uz= \sum_{i=1}^k z_i u_j \in \CX$ follows
  \begin{align*}
    \|\vek Uz\|^2_{\CX} &= \sum_{i=1}^k\sum_{j=1}^k z_iz_j\langle u_i, u_j\rangle_{\CX}\\
    &= \langle z, (\vek U,\vek U)z \rangle_{\R^k}\le \|(\vek U,\vek U )z\|_{\R^k}\|z\|_{\R^k}\\
    &\le\|(\vek U , \vek U )\|_F\|z\|^2_{\R^k}.
  \end{align*}
\end{proof}
From Lemma \ref{Unorm} follows in particular
\be
\|\vek U\|_{\R^k\to \CX} \le \sqrt{\|(\vek U,\vek U)\|_F}.
\ee
We will further on also need
\blem
For $y\in\CY$ and $\vek C$ as above holds
\begin{equation}\label{yC_est}
  \|(y,\vek C)\|_{\R^k} \le \|T\|\|y\|_\CY\sqrt{\sum_{l=1}^k\|u_l\|^2}
\end{equation}
\elem
\begin{proof}
  With \eqref{vekmat} yields
  \begin{align}
    \|(y,\vek C)\|_{\R^k}^2&= \sum_{i=1}^k |(y,c_i)_\CY|^2 \le \sum_{i=1}^k\|y\|^2_\CY\|c_i\|_\CY^2\\
    &=\|y\|^2_\CY\sum_{i=1}^k\|Tu_i\|^2\le \|T\|^2_{_{\CX\to\CY}}\|\|y\|^2_{\CY}\sum_{i=1}^k\|u_i\|_{\CX}^2
  \end{align}
\end{proof}
\section{Projected Regularization Method}\label{section.basic-aug-proj}
Let $\CU\subset\CX$ be a $k$-dimensional subspace, and $\CC=T\,\CU\subset\CY$ be its image under the action of $T$. Here it is important that we assume
that the restricted operator $T|_{\CU}$ is invertible, so as to guarantee that $\dim\CU=\dim\CC=k$ holds.  This allows us express clearly the structure of projectors
onto these spaces.
\bprop\label{prop.grammian-invertible}
Assume that the vectors $u_i\in \vek U$ are linear independent, and that
\begin{equation}
  T_{|_{\CU}}: \CU \to Y
\end{equation}
is invertible on its range $\mathcal{R}(T_{|_{\CU}})$. Then $\prn{\vek C,\vek C}\in \mathbb{R}^{k\times k}$ is continuously invertible.
\eprop
\begin{proof}
According to \eqref{def_subspaces}, $\vek C=\bbmat c_{1} & c_{2} & \cdots & c_{k} \ebmat\in\CY^{k}$ with $c_i=Tu_i$. As $T$ is invertible on $\mathcal{R}(T_{|_{\CU}})$, the vectors in $\vek C$ are linear independent if the vectors in $\vek U$ are linear independent. Moreover,
\be
  \prn{\vek C,\vek C}= \left (\prn{c_i,c_j}\right )_{i,j=1,\dots ,k}
\ee
is the Gramian matrix of the linear independent system contained in $\vek C$ and therefore an invertible matrix. Additionally, as $\prn{\vek C,\vek C}:\mathbb{R}^k\to \mathbb{R}^k$ maps between finite dimensional spaces, its inverse is also bounded.
\end{proof}
In general, the operator $T$ is ill posed and might have a non-trivial null space, i.e., $T$ is not continuously invertible. However, on a finite dimensional subspace invertibility might be achieved easily. \\

Now let $Q$ be the $\CY$-orthogonal projector onto $\CC$ and $P$ {\color{black} be the 
$\prn{\cdot,T^{\ast}T\,\cdot}_{\CX}$-orthogonal} projector onto $\CU$.  The orthogonal projector $Q$ induces a direct sum 
splitting $\CY = \CC \oplus \CC^{\perp_{\prn{\cdot,\cdot}_{\CY}}}$.
 Consider $w\in\CY$.  We represent the action of $Q$ as $Qw = \sum_{i=1}^{k}\prn{w,c_{i}}_{\CY}c_{i}$.
 Similarly, $P$ induces a direct sum splitting of $\CX$ with respect to $\prn{\cdot,T^{\ast}T\,\cdot}_{\CX}$,
 namely $\CX = \CU \oplus \CU^{\perp_{\prn{\cdot,T^{\ast}T\,\cdot}_{\CX}}}$ where $\CU^{\perp_{\prn{\cdot,T^{\ast}T\,\cdot}_{\CX}}}$ denotes
 the orthogonal complement of $\CU$ with respect to the inner product $\prn{\cdot,T^{\ast}T\,\cdot}_{\CX}$.
 We represent the action of $P$ as $Pv = \sum_{i=1}^{k}\prn{v,T^{\ast} T\,u_{i}}_{\CX}u_{i}$.
We also can construct these projectors explicitly.

\bprop\label{eqn.proj-structure}
Let $\vek U\in\CX^{k}$ be a basis of $\CU$ and
$\vek C=T\,\vek U\in\CY^{k}$ be a basis of $\CC$.  Then we have
\be\nn
	P\cdot  = \vek U\prn{\vek U,T^{\ast}T\,\vek U}_{\CX}^{-1}\prn{\cdot,T^{\ast}T\,\vek U}_{\CX}\mand Q\cdot = \vek C\prn{\vek C,\vek C}_{\CY}^{-1}\prn{\cdot,\vek C}_{\CY}.
\ee
\eprop
\begin{proof}
	This is a standard result, particularly in the setting where $\CX=\R^{m}$ and $\CY=\R^{n}$, but it is helpful to see why this is true in this more general setting, as well.
	From \Cref{prop.grammian-invertible}, we know that $\prn{\vek C,\vek C}=\prn{T\vek U,T\vek U} = \prn{\vek U,T^{\ast}T\vek U}$ is continuously invertible.
	Let $\curl{\hat{u}_{k+1},\hat{u}_{k+2},\ldots, \hat{u}_{\ell},\ldots}$ be a basis for $\CU^{\perp_{\prn{\cdot,T^{\ast}T\,\cdot}_{\CX}}}$.  Then for any
	$v\in\CX$, we can write
	\begin{equation}\nonumber
		v = \sum_{i=1}^{k}\prn{v,T^{\ast} T\,u_{i}}_{\CX}u_{i} + \sum_{i=k+1}^{\infty}\prn{v,T^{\ast} T\,\hat{u}_{i}}_{\CX}\hat{u}_{i}.
	\end{equation}
	Since by definition $\prn{\hat{u_{i}},T^{\ast}T\,\vek U}_{\CX}=0$ for all $i>k$, we can write
	\be\label{eqn.proj-scaling-term}
		\vek U\prn{\vek U,T^{\ast}T\,\vek U}_{\CX}^{-1}\prn{v,T^{\ast}T\,\vek U}_{\CX} = \sum_{i=1}^{k}\prn{v,T^{\ast} T\,u_{i}}_{\CX}\vek U\prn{\vek U,T^{\ast}T\,\vek U}_{\CX}^{-1}\prn{u_{i},T^{\ast}T\,\vek U}_{\CX}.
	\ee
	Observe that, according to the definitions in \cref{vekmat} and \cref{matmat}, the vector $\prn{u_{i},T^{\ast}T\,\vek U}_{\CX}$ is the $i$th column of the matrix
	$\prn{\vek U,T^{\ast}T\,\vek U}_{\CX}$.  Thus,
	\be\nn
		\prn{\vek U,T^{\ast}T\,\vek U}_{\CX}^{-1}\prn{u_{i},T^{\ast}T\,\vek U}_{\CX} = \vek e_{i}\in\R^{k}
	\ee
	the $i$th
	Cartesian basis vector.  Thus, we have
	\begin{equation}
		\vek U\prn{\vek U,T^{\ast}T\,\vek U}_{\CX}^{-1}\prn{u_{i},T^{\ast}T\,\vek U}_{\CX} = \vek U\vek e_{i} = u_{i}
	\end{equation}
	and thus according to \cref{eqn.proj-scaling-term}
	\begin{equation}
	\vek U\prn{\vek U,T^{\ast}T\,\vek U}_{\CX}^{-1}\prn{v,T^{\ast}T\,\vek U}_{\CX} = \sum_{i=1}^{k}\prn{v,T^{\ast} T\,u_{i}}_{\CX}u_{i}.
	\end{equation}
	As $v$ was chosen arbitrarily, this is true for any element of $\CX$ meaning the action of $\vek U\prn{\vek U,T^{\ast}T\,\vek U}_{\CX}^{-1}\prn{\cdot,T^{\ast}T\,\vek U}_{\CX}$ is the action of $P$.  This completes the proof for $P$.  The same line of reasoning yields the proof for $Q$.
\end{proof}
 We assume that with exact data that \cref{illposedeq} is consistent. Using the projector $P$, we decompose
the exact true solution,
\be\nn
	x^\dagger = Px^\dagger + \prn{I_{\CX} - P}x^\dagger.
\ee
With exact data, $Px^\dagger$ would be exactly computable, since we can rewrite
\be\label{Pxdef}
	\xpdag:=Px^\dagger = \vek U\prn{\vek C,\vek C}_{\CY}^{-1}\prn{Tx,\vek C}_{\CY} = \vek U\prn{\vek C,\vek C}_{\CY}^{-1}\prn{y,\vek C}_{\CY}.
\ee
However, we do no have exact data.  We have disturbed data $y^{\delta}$ such that $\nm{y-y^{\delta}}_{\CY} \le \delta$.  Rewriting
\be\label{PXrewritten}
Px = \vek U\prn{\vek C,\vek C}_{\CY}^{-1}\prn{y-y^{\delta},\vek C}_{\CY} + \vek U\prn{\vek C,\vek C}_{\CY}^{-1}\prn{y^{\delta},\vek C}_{\CY}
\ee
we can define
\begin{equation}\label{xdelta}
  \xpd := \vek U\prn{\vek C,\vek C}_{\CY}^{-1}\prn{y^{\delta},\vek C}_{\CY},
\end{equation}
which can be computed {\color{black} inexpensively to high precision} and may be used as an approximation $\xpd\approx Px^\dagger$. We obtain
\bprop\label{prop.approxdiff}
The approximation error between $\xpdag$ and $\xpd$ is given by
\be\label{approxdiff}
	\xpdag - \xpd = \vek U\prn{\vek C,\vek C}_{\CY}^{-1}\prn{y-y^{\delta},\vek C}_{\CY},
\ee
which can be estimated as
\be\label{eqn.init-proj-err}
\nm{\xpdag - \xpd }\le	\left (\sqrt{\|(\vek U,\vek U)\|_F \sum_{l=1}^k\|u_l\|^2}\nm{\prn{\vek C,\vek C}^{-1}}_{F}\right )\delta .
\ee
\eprop
\begin{proof}
Using \eqref{Pxdef}-\eqref{xdelta} gives \eqref{approxdiff}, which can be estimated as
\begin{align}
  \nm{\xpdag - \xpd}_{\CX} &\le  \nm{\vek U\prn{\vek C,\vek C}_{\CY}^{-1}\prn{y-y^{\delta},\vek C}_{\CY}}\\
  &\le \nm{\vek U}_{\R^k\to \CX}\nm{\prn{\vek C,\vek C}^{-1}}_{\R^k\to\R^k}\nm{\prn{y-y^{\delta},\vek C}_{\CY}}_{\R^k}\\
  &\stackrel{\eqref{Unorm},\eqref{yC_est}}{\le} \sqrt{\|(\vek U,\vek U)\|_F}\nm{\prn{\vek C,\vek C}^{-1}}_{\R^k\to\R^k} \|T\|\|y-y^\delta\|_\CY\sqrt{\sum_{l=1}^k\|u_l\|^2}\\
  &= \left (\sqrt{\|(\vek U,\vek U)\|_F \sum_{l=1}^k\|u_l\|^2}\nm{\prn{\vek C,\vek C}^{-1}}_{F}\right )\delta
\end{align}
\end{proof}

\brem
Often with such recycling methods, one takes either $\CU$ or $\CC$ to be represented with an orthonormal basis, for the purposes of convenience with
regard to analysis and implementation.  We assume here that $\CC$ is represented by an orthonormal basis.
If we scale $\vek U$ such that $\vek C=T\,\vek U$ is an orthonormal system spanning $\CC$, we can write \cref{eqn.init-proj-err} more compactly as
\be\nn
	\nm{Px - \xpd}_{\CX} \leq \nm{\vek U}_{\CX} \delta = \CO\prn{\delta}.
\ee
Furthermore, $P$ and $Q$ can be expressed more compactly,
\be\nn
	P\cdot  = \vek U\prn{\cdot,T^{\ast}T\,\vek U}_{\CX}\mand Q\cdot = \vek C\prn{\cdot,\vek C}_{\CY},
\ee
due to the fact that $\prn{\vek U,T^{\ast}T\,\vek U}_{\CX}=\prn{\vek C,\vek C}_{\CY} = \vek I$ when $\vek C$ is an orthonormal system.
\erem

It remains to find a suitable approximation for $ \prn{I_{\CX} - P}x^\dagger$. We have the following result:
\bprop
Under the assumption
\begin{equation}\label{eq:ProjAssumpt}
  QT(I-P)x^\dagger = 0
\end{equation}
holds
\begin{equation}\label{eq:ProjEqCompl}
  y-y_p =(I-Q)T(I-P)x^\dagger ,
\end{equation}
where
\begin{equation}
  y_p:= T\xpdag .
\end{equation}
\eprop
\begin{proof}
  We have
  \begin{equation}\label{eq:Txddecomp}
  y=  Tx^\dagger = QTPx^\dagger +QT(I-P)x^\dagger + (I-Q)TPx^\dagger +(I-Q)T(I-P)x^\dagger .
\end{equation}
  As $\xpd=Px^\dagger\in\CU$ we have $y_p=T\xpd\in T\CU$ and therefore
  \begin{equation}\label{eq:QTPx}
    y_p=T\xpd = QTPx^\dagger .
  \end{equation}
  Further on, as $TPx^\dagger\in T\CU$ and $(I-Q)$ projects onto the complement of $T\CU$ we obtain
  \begin{equation}\label{eq:IQTP}
    (I-Q)TPx^\dagger =0.
  \end{equation}
  Inserting  \eqref{eq:QTPx}, \eqref{eq:IQTP} into \eqref{eq:Txddecomp} yields
  \begin{equation}
    y-y_p = QT(I-P)x^\dagger + (I-Q)T(I-P)x^\dagger ,
  \end{equation}
  and Assumption \eqref{eq:ProjAssumpt} yields \eqref{eq:ProjEqCompl}
\end{proof}
Assumption \eqref{eq:ProjAssumpt} transforms to
\begin{eqnarray}\nonumber
  QT(I-P)x^\dagger = 0 &\Longleftrightarrow & \langle z, QT(I-P)x^\dagger \rangle = 0 \hspace{1cm}\forall z\in \CY \\
  &\Longleftrightarrow & \langle T^\ast Qz, (I-P)x^\dagger \rangle = 0\hspace{1cm}\forall z\in \CY \label{eq:orthosol}
\end{eqnarray}
We arrive at
\bprop
  If
  \begin{equation}\label{eq:PerpCond}
    T^\ast \CC \perp (I-P)x^\dagger ,
  \end{equation}
  then \eqref{eq:ProjAssumpt} holds. Particularly, \eqref{eq:PerpCond} is fulfilled  if
  \begin{equation}\label{eq:TstarTinclusion}
    T^\ast T \CU \subset \CU
  \end{equation}
  holds.
\eprop
\begin{proof}
  As $Q\CY=\CC$, eq. \eqref{eq:orthosol} transforms to \eqref{eq:PerpCond}. As $(I-P)x^\dagger \in \CU^\perp$ holds always if we additionally assume that
  $\CU$ is chosen s.t. $T^\ast \CC \perp \CU^\perp$,
   and with $\CC = T\CU$ this transforms to $T^\ast T \CU\perp \CU^\perp$ or, equivalently, \eqref{eq:TstarTinclusion}.
  \end{proof}
\begin{remark}
Obviously, \eqref{eq:TstarTinclusion} is much stronger than \eqref{eq:PerpCond}. However, there are easy examples where \eqref{eq:PerpCond} always holds:
\begin{enumerate}
  \item If $\CU = span \{u_i: i=1,\dots N \mbox{ and } T^\ast T u_i = \sigma_i u_i\}$ then clearly $T^\ast T \CU\subset \CU$.
  \item If $T$ is an {\it unitary}  operator, i.e., $\langle Tx_1,Tx_2\rangle=\langle x_1,x_2\rangle$ holds for all $x_1,x_2$, then we have for $u\in \CU$ and $u^\perp\in \CU^\perp$
  \[\langle T^\ast T u, u^\perp\rangle= \langle T u, Tu^\perp\rangle=\langle u, u^\perp\rangle=0 ,\]
  i.e.,   \eqref{eq:TstarTinclusion}.
\end{enumerate}
\end{remark}
Furthermore, the way we construct our projectors ensures that the assumption holds.
\blem
	The assumption \cref{eq:ProjAssumpt} holds for any pair of projectors $\prn{P,Q}$ satisfying \cref{eqn.projector-relationship}.
\elem
\bproof
	This follows from the fact that the product of complementary projectors is the zero operator, since from \cref{eqn.projector-relationship}, we have
	\be\nn
		 QT(I-P)x^\dagger = Q(I-Q)Tx^\dagger =  0.
	\ee
\eproof
With the above results we can now compute $(I-P)x^\dagger$ as a solution of equation \eqref{eq:ProjAssumpt}. However, as always in Inverse Problems, we might not have the exact right hand side $y-y_p$ but some noisy version $y^\delta -y_p^\delta$, where we define $y_p^\delta := T\xpd$. The data error can be approximated as follows:

\bprop
  The data error in the left hand side of \eqref{eq:ProjAssumpt} can be estimated as
  \begin{equation}\label{eq:proj_dataerror}
\|(y-y_p) - (y^\delta -y_p^\delta)\|\le \kappa_{\CU}\cdot\delta
  \end{equation}
  with
  \begin{equation}
    \kappa_{\CU}:= 1+\|T\|\sqrt{\|(\vek U,\vek U)\|_F \sum_{l=1}^k\|u_l\|^2}\nm{\prn{\vek C,\vek C}^{-1}}_{F}
  \end{equation}
\eprop
\begin{proof}
  The proof is straight forward:

  \begin{equation}
    \|(y-y_p) - (y^\delta -y_p^\delta )\|\le  \|y-y^\delta\| + \|y_p-y_p^\delta\|\le\delta +\|T\|\|\xpdag - \xpd\|
  \end{equation}
  and with \eqref{eqn.init-proj-err} follows \eqref{eq:proj_dataerror}
\end{proof}

We are now ready to formulate our solution approach:

\balgte[H]\label{alg.aug}
\caption{Augmented Regularization}
Given:\newline
  $\bullet$ data $y^\delta$ fulfilling $\|y-y^\delta\|\le \delta$\newline
$\bullet$ $\CU\subset \CX$, $\CU$ finite dimensional\newline
  $\bullet$ a regularization method $R_\alpha (y,T)$ for the equation $y=Tx$ with parameter choice rule $\alpha=\alpha (\delta)$\\
Compute $x_p^\delta$ according to \eqref{xdelta} and set $y_p^\delta = Tx_p^\delta$\\
Set $B:=(I-Q)T$\\
Compute $\hat{x}_p^\delta:= R_\alpha (y^\delta-y_p^\delta,B)$, $\alpha=\alpha (\kappa_{\CU}\delta)$ \\
Set $x_{\CU,\alpha}^\delta:=x_p^\delta+\hat{x}_p^\delta$
\ealgte

\bprop
If for the solution $x^\dagger$ of the equation $Tx=y$ and the chosen augmentation subspace $\CU$ eq. \eqref{eq:ProjAssumpt} holds, then the
 augmented regularization with data $y^\delta$ fulfilling $\|y-y^\delta\|\le \delta$ as described in Algorithm \ref{alg.aug} forms a regularization method.
\eprop
\begin{proof}
 As $R_\alpha (y^\delta-y_p^\delta ,B)$ is a regularization for \eqref{eq:ProjEqCompl} we have
 as $\delta\to 0$, and decomposing $x^\dagger = Px^\dagger+(I-P)x^\dagger=x_p^\dagger + \hat{x}_p^\dagger$ yields
   \begin{eqnarray*}
     \|x^\dagger - x_{\CU,\alpha}^\delta\|&\le & \|x_p^\dagger -x_p^\delta\|+\|\hat{x}_p^\dagger -\hat{x}_p^\delta\|\\
     &\stackrel{\eqref{eqn.init-proj-err}}{\le}&\mathcal{O}(\delta)+\|\hat{x}_p^\dagger - R_{\alpha (\kappa_{\CU}\delta)} (y^\delta -y_p^\delta,B)\|\stackrel{\delta \to 0}{\longrightarrow} 0
   \end{eqnarray*}
   which concludes the proof.
\end{proof}
\brem
	We can thus extend any existing regularization $R_\alpha$ to its augmented method and still obtain a regularization method.
\erem

\bex
Let us illustrate this for the {\it Landweber method}, which has been defined in \eqref{Landweber}. If the iteration is stopped by the discrepancy principle, i.e., if the stopping index $k_\ast$ is determined as the first index such that
\begin{equation}
  \|y^\delta -Tx_{k_\ast}^\delta\|\le \tau \delta \hspace{1cm} \tau >1
\end{equation}
holds, then it is a regularization method, see e.g., \cite{demol:1}. For given data $y^\delta$, the augmented Landweber method would read as
\begin{equation}
  x_{k_\ast}^\delta = x_p^\delta + \hat{x}_{k_\ast}^\delta
\end{equation}
where $x_p^\delta$ is defined in \eqref{xdelta} and $\hat{x}_{k_\ast}^\delta $ is the $k_\ast$th iterate of the Landweber iteration applied to the equation
\begin{equation}
(I-Q)Tx= y^\delta- y_p^\delta, \hspace{2cm} y_p^\delta = T x_p^\delta ~.
\end{equation}
Using
\begin{equation}
  y_p^\delta = Tx_p^\delta = T \vek U\prn{\vek C,\vek C}_{\CY}^{-1}\prn{y^{\delta},\vek C}_{\CY} = \vek C\prn{\vek C,\vek C}_{\CY}^{-1}\prn{y^{\delta},\vek C}_{\CY}= Q y^\delta
\end{equation}
the iteration for $\hat{x}_k^\delta$ reads as
\begin{eqnarray}\nonumber
  \hat{x}_{k+1}^\delta &=& \hat{x}_k^\delta+T^\ast (I-Q)\left (y^\delta-  y_p^\delta - (I-Q)T  \hat{x}_k^\delta\right )\\
  &=& \hat{x}_k^\delta+T^\ast (I-Q)\left (y^\delta-T  \hat{x}_k^\delta\right ).
\end{eqnarray}
Usually such an iteration is started with $\hat{x}_0^\delta=0$. Instead, we can incorporate $x_p^\delta$ directly into the iteration by setting $\hat{x}_0^\delta=x_p^\delta$.
\eex
{\color{black}
\subsection{What to recycle}\label{sec.what-to-recycle}
One question we have not addressed thus far concerns what should the subspace $\CU$ encode?
This is a question that hinges very much on the application.  Suffice it to say,  a detailed review of recycling
strategies is beyond the scope of this paper, but we refer the reader to \cite[Section 6]{soodhalter2020survey} for
a more detailed discussion.  Some common choices for augmentation vectors are approximate solutions or
solutions to related problems, approximate or (when available) exact eigenvectors/singular vectors, and vectors satisfying
convergence model optimality properties.  

The use of approximate eigenvectors as a recycling strategy was first suggested for recycling between cycles of GMRES
applied to one linear system to mitigate the effects on convergence speed of restarting \cite{morgan.gmresdr}.  
The algorithmic
choices in that paper necessitated that the approximate eigenvectors be harmonic Ritz vectors.  Ritz-type 
eigenvector approximations with respect to a square matrix $\vek A$ 
are obtained applying a Galerkin or Petrov-Galerkin approximation strategy
to the eigenvalue problem of the form
\begin{equation}\nonumber
	\mbox{select}\qquad \vek v\in\mathcal{V}\qquad \mbox{such that}\qquad \vek A \vek v-\lambda \vek v\perp \mathcal{W},
\end{equation}
for some $\lambda\in\R$.  The spaces $\mathcal{V}$ and $\mathcal{W}$ are usually spaces which have been generated
during the iteration, so that the solution to the Ritz problem reduces to a (generalized) eigenvalue problem.  This strategy
was extended and combined with GCRO-type optimal methods in a way that did not require augmentation by
harmonic Ritz vectors, and this enabled recycling between multiple linear systems \cite{Parks.deSturler.GCRODR.2005}.
This has been further extended to the case of short-recurrence schemes not requiring restart; see, e.g., \cite{dSC.2019,Kilmer.deSturler.tomography.2006,WSG.2007}.  In that case, one 
 stores and updates a running recycle \textit{window} but not to use in the current iteration.  Rather, it is being built to
 be used in subsequent systems.  
 
It should be noted that in the context of discrete ill-posed problems, a poor or noisy choice of 
$\CU$ may induce a projected problem \cref{eqn.proj-prob} which is even more ill-posed.  However, 
it has been demonstrated that there are many good choices for $\CU$; see, e.g., \cite{BR-2.2007,BR.2007,DGH.2014}.  
More recently,
a recycled LSQR algorithm has
been proposed where approximate eigenvectors are used as well as other data which become available in real-time
\cite{dSC.2019}.  As LSQR is a short-recurrence method, the authors employ a windowed recycling strategy, where the
recycled subspace is updated every $p$ iterations so that only a window of vectors generated by the short recurrence
iteration must be stored.  They propose a
few different recycling approaches appropriate for ill-posed problems: truncated singular vectors, solution- or sparsity-oriented compression, and reduced basis decomposition.

In the next section, we use the recycling method framework from \cite{soodhalter2020survey} to develop augmented
versions of steepest descent and Landweber methods.  In \Cref{sec.numerics}, we combine this method with recycling
strategies appropriate to each individual application, building on strategies from the literature but adjusted to fit the
particular applications.

}

\section{Augmented gradient descent methods}
We now delve further into the case that we apply a gradient descent-type method (i.e., Landweber) to the projected problem, in order
to produce a practical implementation of this algorithm.  In particular, we first explore augmenting the steepest 
descent method, which can be understood as a Landweber-type method with a step-length that is dynamically 
chosen to minimize the residual norm.  We have shown in \Cref{eqn.steepest-descent-characterization} that this is
a projection method and thus can be augmented in the above-discussed framework. 
Augmented methods such as this rely on certain identities to achieve algorithmic advantages.
\blem\label{lem.proj-resid-equiv}
	The residual produced by the augmented iterative method with approximation $\hat{x}_k^\delta = \vek U\prn{\vek C,\vek C}_{\CY}^{-1}\prn{y^{\delta},\vek C}_{\CY} + \prn{I_{\CX}-P}\hat{t}_k^\delta$ defined as in \cref{eqn.x-approx-decomp}
	where $\hat{t}_k^\delta$ is the approximate solution of \cref{eqn.proj-prob} coincides with the
	residual associated to $\hat{t}_k^\delta$ for the projected problem \cref{eqn.proj-prob}, i.e.,
	\be\nn
		y^{\delta} - T\hat{x}_k^\delta = \prn{I_{\CY} - Q}\prn{y^{\delta} - T\hat{t}_k^\delta} \perp\CC.
	\ee
\elem
\bproof
	This result has been shown in a number of papers for finite-dimensional, well-posed problems; 
	see, \cite{dSKS.2018} and \cite{deSturler.GCRO.1996,Gaul.2014-phd,GGL.2013,Gutknecht.AugBiCG.2014}.  We show it here for completeness by 
	direct calculation, by first writing
	\begin{align*}
			y^{\delta} - T\hat{x}_j^\delta &= y^{\delta} - T\prn{\vek U\prn{\vek C,\vek C}_{\CY}^{-1}\prn{y^{\delta},\vek C}_{\CY} + \prn{I_{\CX}-P}\hat{t}_j^\delta}\\
			y^{\delta} - T\hat{x}_j^\delta &= y^{\delta} - \prn{\vek C\prn{\vek C,\vek C}_{\CY}^{-1}\prn{y^{\delta},\vek C}_{\CY} + T\prn{I_{\CX}-P}\hat{t}_j^\delta}.
	\end{align*}	
	Then we observe that it follows from \Cref{eqn.proj-structure} and \cref{eqn.projector-relationship} that this is equivalent to
	\begin{align*}
			y^{\delta} - T\hat{x}_j^\delta &= y^{\delta} - \prn{Qy^{\delta} + \prn{I_{\CX}-Q}T\hat{t}_j^\delta}\\
			&= \prn{I_{\CX}-Q}y^{\delta} -  \prn{I_{\CX}-Q}T\hat{t}_j^\delta.
	\end{align*}	
	It follows directly that $y^{\delta} - T\hat{x}_j^\delta\perp\CC$.
\eproof
Note that this proof holds in the presence or absence of noise and whether or not the problem being treated is consistent.  In the case that
the right-hand side is noise-free and the problem is consistent, we have from \cref{PXrewritten} that 
 \be\label{eqn.MPsol-proj}
 	\vek U\prn{\vek C,\vek C}_{\CY}^{-1}\prn{y,\vek C}_{\CY} = Px^{\dagger} = x_p^\dagger,
 \ee 
 meaning it 
is the best approximation of $x^{\dagger}$ in $\CU$ with respect
to the $T^{\ast}T$-norm.  In the presence of noise, we previously characterized and bounded the error introduced into \cref{eqn.MPsol-proj}
in \Cref{prop.approxdiff}.

It directly follows that the construction of the gradient descent direction for Landweber applied to \cref{eqn.proj-prob} can be simplified.
\bcor
	The gradient descent direction for Landweber applied to \cref{eqn.proj-prob} admits the simplification,
	\be\nn
		T^{\ast}\prn{I_{\CY} - Q}\prn{y^{\delta} - T\hat{t}^{\delta}_{j}} = T^{\ast}\prn{y^{\delta} - T\hat{x}^{\delta}_{j}};
	\ee
	i.e., the descent direction is the same as that for the full problem \cref{illposedeq}.
\ecor

With these simplifications, we now show that indeed the above-described fully augmented Landweber construction is indeed updating the
approximation to the solution to the full problem \cref{illposedeq} in the optimal direction at each iteration. Furthermore, if when applying Landweber
to \cref{eqn.proj-prob}, we choose the step length dynamically such that we are applying the method of \emph{steepest descent} to \cref{eqn.proj-prob},
this will be shown in the following theorem to be equivalent to minimizing the residual over the sum of subspaces 
$\CS_{j+1} = \CU  + \Span\curl{T^{\ast}\prn{y^{\delta} - Tx_{j}}}$. 

As this should be a steepest descent-type method which is a $\CY$-norm minimization method, the correct constraint space is the image of the correction space under the action of the operator, namely
$\CL_{j+1} = T\CU + \Span\curl{TT^{\ast}\prn{y^{\delta} - Tx_{j}}}$.  
At step $j+1$, we have $x_{j+1} = x_{j} + \vek U\vek u_{j+1} + \alpha_{j+1}T^{\ast}\prn{y^{\delta} - Tx_{j}}$ where $\vek u_{j+1}\in\R^{k}$
and $\alpha_{j+1}\in\R$ are determined by the minimization constraint.
The following is an adaption of a more general result, proven in \cite{dSKS.2018}. 
\bthm
	The $(j+1)$st approximate solution to the augmented steepest descent problem minimization,
	\begin{align*}
		&\mselect \vek U\vek u_{j+1} + \alpha_{j+1}v_{j}\in\CU + \Span\curl{v_{j}}\mwhere v_{j} = T^{\ast}\prn{y^{\delta} - Tx_{j}} \\&\msuchthat \prn{\vek u_{j+1}, \alpha_{j+1}} = \argmin{\vek u\in\R^{k}\atop\alpha\in\R}\nm{y - T\prn{x_{j} + \alpha v_{j} + \vek U\vek u}}_{\CY}
	\end{align*}
	satisfies
	\begin{enumerate}
		\item $\vek u_{j+1} = \vek u^{(1)}_{j+1} +\vek u^{(2)}_{j+1} $ with 
		\begin{align}
			\vek U\vek  u^{(1)}_{j+1} &= \begin{cases}
																\argmin{x\in\CU}\nm{e^{\dagger}_{0} - x}_{\prn{T^{\ast}T\cdot,\cdot}_{\CX}} - u_{y} &\qquad \mif j=0\\
																\vek 0 & \qquad \motherwise 
															\end{cases}\label{eqn.error-proj}\\
	\vek U\vek u^{(2)}_{j+1} &= - \argmin{x\in\CU}\nm{x - \alpha_{j+1} v_{j+1}}_{\prn{T^{\ast}T\cdot,\cdot}_{\CX}}\label{eqn.v-proj}
		\end{align}
		where $u_{y} = \vek U\prn{\vek C,\vek C}_{\CY}^{-1}\prn{y-y^{\delta},\vek C}_{\CY}$,
		\item and $\alpha_{j+1}$ satisfies
		\be\label{eqn.proj-min-prob}
			\alpha_{j+1}v_{j+1} = \argmin{x\in\Span\curl{v_{j}}}\nm{\prn{I_{\CY} - Q}r_{j} - \prn{I_{\CY} - Q}Tx}_{\CY},
		\ee
	\end{enumerate}
	i.e., $\alpha_{j+1}$ is obtained via applying an iteration of steepest descent to the projected problem, 
	where $v_{j+1} = T^{\ast}r_{j}$ and $e^{\dagger}_{0} = x^{\dagger} - x_{0}$.
\ethm
\bproof
As this is a minimization problem, we prove the result by studying the structure of the associated orthogonal projector.
We can rewrite this minimization as 
\be\nn
 \prn{\vek u_{j+1}, \alpha_{j+1}} = \argmin{\vek u\in\R^{k}\atop\alpha\in\R}\nm{r_{j} - T\prn{ \alpha v_{j} + \vek U\vek u}}_{\CY}.
\ee
Thus, we seek the pair $\prn{\vek u_{j+1}, \alpha_{j+1}}$ giving the $\CY$-norm best approximation of $r_{j}$ in the space
$T\CS_{j+1}$.  This is equivalent to computing the $\CY$-orthogonal projection of $r_{j}$ into that space, which we denote $Q_{T\CS_{j+1}}$.  
Recalling that $\vek C = T\vek U$,
from \Cref{eqn.proj-structure}, we can express this as
\begin{align*}
	T\prn{ \alpha v_{j} + \vek U\vek u} &= Q_{T\CS_{j+1}}r_{j} \\
	&= \footnotesize\bbmat \vek C & Tv_{j+1} \ebmat\prn{\bbmat \vek C & Tv_{j+1} \ebmat,\bbmat \vek C & Tv_{j+1} \ebmat}_{\CY}^{-1}\prn{r_{j},\bbmat \vek C & Tv_{j+1} \ebmat}_{\CY}\normalsize,
\end{align*}
where $\bbmat \vek C & Tv_{j+1} \ebmat\in\CL\prn{\R^{(k+1)},\CY}$.  This means that $\prn{\vek u_{j+1}, \alpha_{j+1}}$ form the solution of
\be\nn
	\bbmat 
			\prn{\vek C,\vek C}_{\CY} & \prn{T v_{j+1},\vek C}_{\CY} \\ 
			\prn{\vek C,T v_{j+1}}_{\CY} & \prn{T v_{j+1},T v_{j+1}}_{\CY}	
	  \ebmat 
	  \bbmat 
	  		\vek u_{j+1} \\ \alpha_{j+1} 
	  \ebmat 
	  = 
	  \bbmat  
	  		\prn{r_{j},\vek C}_{\CY}\\
	  		\prn{r_{j},Tv_{j+1}}_{\CY}
	  \ebmat
\ee
To prove the statement of the theorem, we eliminate $\vek u_{j+1}$ from the second block of equations (following \cite{parks2016block}) yielding
the two sets of equations,
\begin{align*}
	\prn{\vek C,\vek C}_{\CY}\vek u_{j+1} + \alpha_{j+1} \prn{T v_{j+1},\vek C}_{\CY} &= \prn{r_{j},\vek C}_{\CY}\\
	\alpha_{j+1}\left[ \prn{T v_{j+1},T v_{j+1}}_{\CY} -\prn{\vek C,T v_{j+1}}_{\CY} \prn{\vek C,\vek C}_{\CY}^{-1}\prn{T v_{j+1},\vek C}_{\CY}\right]
	 &= \prn{r_{j},Tv_{j+1}}_{\CY} 
	\\ -\prn{\vek C,T v_{j+1}}_{\CY}& \prn{\vek C,\vek C}_{\CY}^{-1}\prn{r_{j},\vek C}_{\CY}.
\end{align*}
We can consolidate the second set of equations, yielding
\begin{align*}
	\alpha_{j+1} \prn{T v_{j+1}-\vek C\prn{\vek C,\vek C}_{\CY}^{-1}\prn{T v_{j+1},\vek C}_{\CY},T v_{j+1}}_{\CY} 
	 = 
	\\   \prn{r_{j}-\vek C\prn{\vek C,\vek C}_{\CY}^{-1}\prn{r_{j},\vek C}_{\CY},T v_{j+1}}_{\CY}.
\end{align*}
Recalling the definition of the projector $Q$, we substitute and then take advantage of the self-adjointness of orthogonal
projectors,
\begin{align*}
	\alpha_{j+1} \prn{\prn{I_{\CY} - Q}T v_{j+1},T v_{j+1}}_{\CY}  &=  \prn{\prn{I_{\CY} - Q}r_{j},T v_{j+1}}_{\CY}\\
	\alpha_{j+1} \prn{\prn{I_{\CY} - Q}^{2}T v_{j+1},T v_{j+1}}_{\CY}  &=  \prn{\prn{I_{\CY} - Q}^{2}r_{j},T v_{j+1}}_{\CY}\\
	\alpha_{j+1} \prn{\prn{I_{\CY} - Q}T v_{j+1},\prn{I_{\CY} - Q}T v_{j+1}}_{\CY}  &=  \prn{\prn{I_{\CY} - Q}r_{j},\prn{I_{\CY} - Q}T v_{j+1}}_{\CY}.
\end{align*}
Thus, the minimizer $\alpha_{j+1}$ satisfies
\be\nn
	\bigprn{\prn{I_{\CY} - Q}r_{j} - \prn{I_{\CY} - Q}T \prn{\alpha_{j+1}v_{j+1}},\prn{I_{\CY} - Q}T v_{j+1}}_{\CY} = 0.
\ee
As $v_{j+1} = T^{\ast}r_{j}$, then we see that $\alpha_{j+1}$ is the solution of
\small\begin{align*}
	&\mselect \alpha_{j+1} v_{j+1}\in \Span\curl{T^{\ast}\prn{I_{\CY} - Q}r_{j}}\\ 
	&\msuchthat \prn{I_{\CY} - Q}r_{j} - \prn{I_{\CY} - Q}T \prn{\alpha_{j+1}v_{j+1}}\perp \Span\curl{\prn{I_{\CY} - Q}TT^{\ast}\prn{I_{\CY} - Q}r_{j}}.
\end{align*}\normalsize
From \Cref{eqn.steepest-descent-characterization}, we then see that this is the residual constraint formulation of the $(j+1)$st step of steepest descent
applied to \cref{eqn.proj-min-prob}.
We can then solve for $\vek u_{j+1}$,
\be\nn
	\vek u_{j+1} = \prn{\vek C,\vek C}_{\CY}^{-1}\brac{\prn{r_{j},\vek C}_{\CY} -  \alpha_{j+1} \prn{T v_{j+1},\vek C}_{\CY}}.
\ee
Expanding to obtain the contribution from $\CU$ yields
\begin{align*}
	\vek U \vek u_{j+1} = \vek U\prn{\vek C,\vek C}_{\CY}^{-1}\prn{r_{j},\vek C}_{\CY} -  \alpha_{j+1}  \vek U\prn{\vek C,\vek C}_{\CY}^{-1}\prn{T v_{j+1},\vek C}_{\CY},
\end{align*}
and we observe that second term in this sum is \cref{eqn.v-proj}.  For the first term, we must consider two cases.  If $j=0$, then
$x_{0} = 0$, and $r_{0}^{\delta} = y^{\delta}$.  Thus, we can write
\be\nn
	\vek U\prn{\vek C,\vek C}_{\CY}^{-1}\prn{r_{0},\vek C}_{\CY} = \vek U\prn{\vek C,\vek C}_{\CY}^{-1}\prn{y^{\delta},\vek C}_{\CY} = x_{p}^{\delta}.
\ee
From \Cref{prop.approxdiff}, we have that $x_{p}^{\delta}  = x_{p}^{\dagger} - U\prn{\vek C,\vek C}_{\CY}^{-1}\prn{y-y^{\delta},\vek C}_{\CY}$.
It is established in \cref{eqn.MPsol-proj} that $x_{p}^{\dagger} $ is the minimizer in \cref{eqn.error-proj}.
Lastly, if $j>0$, \Cref{lem.proj-resid-equiv} tells us that $r_{j}\perp\CC$; thus $\prn{r_{j},\vek C}_{\CY}=\vek 0$.
\eproof

\bcor
	The full augmented steepest descent approximation can be represented as 
	\be\nn
		x_{j+1} 
		= 
		\begin{cases}
			x_{p}^{\delta} + \alpha_{1} T^{\ast}y^{\delta}  - \alpha_{i} P\prn{ T^{\ast}y^{\delta} } &\mif j=0 \\
			x_{j} + \alpha_{i} T^{\ast}r_{j}  - \alpha_{i} P\prn{ T^{\ast}r_{j}}& \motherwise
		\end{cases}
	\ee
	with optimal step length
	\be\nn
		\alpha_{i} = \dfrac{\nm{T^{\ast}r_{i}}_{\CY}^{2}}{\nm{\prn{I_{\CY} - Q}TT^{\ast}r_{i}}_{\CY}^{2}}.
	\ee
\ecor

We note that $x_{p}^{\delta}\in\CU$ can be computed one time initially, as it comes only from initial data.  
We encapsulate all of this in the case of augmented steepest descent with dynamically determined optimal step length 
as Algorithm \ref{alg.aug-steepest-descent}. 

\textbf{Note}:  the reader may notice that Line 3 of Algorithm \ref{alg.aug-steepest-descent}
asks that the user compute a {\tt QR}-factorization of $T\vek U$, which essentially has a finite number of 
elements from an infinite-dimensional space $\CY$ as its ``columns''.  By {\tt QR}-factorization in this context, we mean 
performing steps of the Gram-Schmidt process and storing the orthogonalization and normalization coefficients in an
upper-triangular matrix $\vek R$.  Right-multiplication of $\vek C\in\CL\prn{\R^{k},\CY}$ by a matrix should be understood
as applying the definition in \cref{UC-Def} for each column of $\vek R$, thereby recovering the relationship
between $T\vek U$ and $\vek C$ arising from the Gram-Schmidt process. 
This is a well-defined procedure and process which is discussed, e.g., in \cite[Lecture 7]{trefethen.bau}.

\balgte[H]\label{alg.aug-steepest-descent}
\caption{Augmented steepest descent for the normal equations}
Given: $\vek U\in\CL\prn{\R^{k},\CX}$ representing $\CU$\\
Set $r_{0} = y - Tx_{0}$\\
Compute {\tt QR}-factorization $T\vek U = \vek C\vek R$\\
$\vek U\leftarrow \vek U\vek R^{-1}$\\
$\vek u^{(1)} = \prn{r_{0},\vek C}_{\CY}$\\
$x \leftarrow x_{0} + \vek U\vek u^{(1)}  $\\
$r \leftarrow r_{0} - \vek C\vek u^{(1)}  $\\
\While{STOPPING-CRITERIA}{
	$\alpha_{i} = \frac{\nm{T^{\ast}r_{i}}_{\CY}^{2}}{\nm{\prn{I_{\CY} - Q}TT^{\ast}r_{i}}_{\CY}^{2}}$\\
	$\widehat{\vek w}\leftarrow = \prn{TT^{\ast}r,\vek C}_{\CY}$\\
	$x\leftarrow x + \alpha_{i}T^{\ast}r - \alpha_{i}\vek U\widehat{\vek w}$\\
	$r\leftarrow r - \alpha_{i}TT^{\ast}r + \alpha_{i}\vek C\widehat{\vek w}$
}
\ealgte

Furthermore, it is then clear that if we no longer choose $\alpha_{i}$ dynamically according to a minimization criteria and instead fix $\alpha_{i}\equiv\alpha$
that we can propose an augmented Landweber method, which we encapsulate as Algorithm \ref{alg.aug-Landweber}.

\balgte[H]\label{alg.aug-Landweber}
\caption{Augmented Landweber for the normal equations}
Given: $\vek U\in\CL\prn{\R^{k},\CX}$ representing $\CU$,\, $\alpha>0$\\
Set $r_{0} = y - Tx_{0}$\\
Compute {\tt QR}-factorization $T\vek U = \vek C\vek R$\\
$\vek U\leftarrow \vek U\vek R^{-1}$\\
$\vek u^{(1)} = \prn{r_{0},\vek C}_{\CY}$\\
$x \leftarrow x_{0} + \vek U\vek u^{(1)}  $\\
$r \leftarrow r_{0} - \vek C\vek u^{(1)}  $\\
\While{STOPPING-CRITERIA}{
	$\widehat{\vek w}\leftarrow = \prn{TT^{\ast}r,\vek C}_{\CY}$\\
	$x\leftarrow x + \alpha T^{\ast}r - \alpha\vek U\widehat{\vek w}$\\
	$r\leftarrow r - \alpha TT^{\ast}r + \alpha\vek C\widehat{\vek w}$
}
\ealgte
~

\section{Numerical experiments}\label{sec.numerics}
We present here experiments from two subproblems from the field of adaptive optics for large ground-based telescopes as well as for a toy blurring problem
presented as a proof-of-concept for an alternative use of these methods.


\subsection{Gaussian blurring model}%

\begin{figure}[H]\label{fig.augSteepestDescent-resid}
	\centering
	\includegraphics[scale=.18]{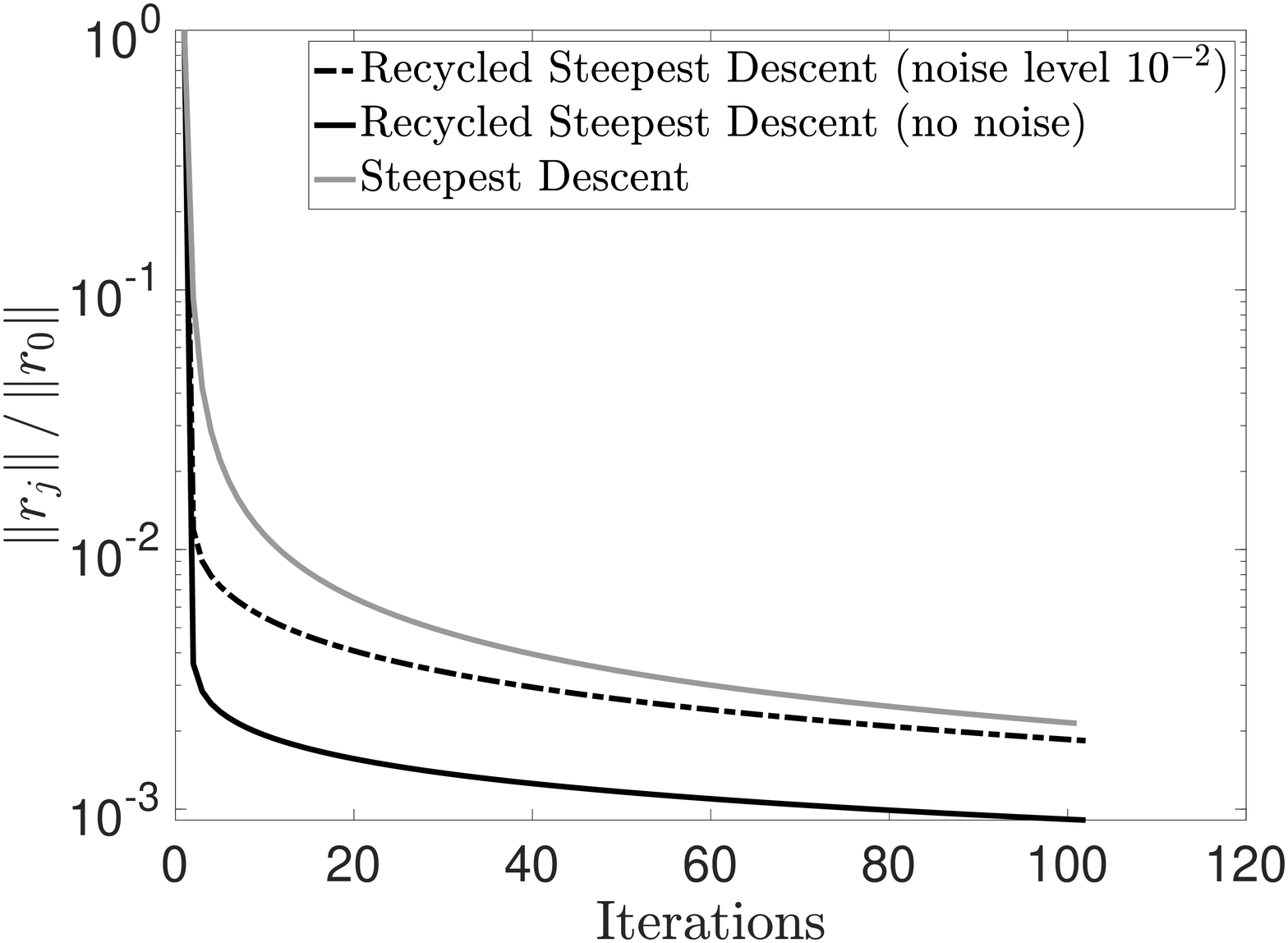}
	\includegraphics[scale=.18]{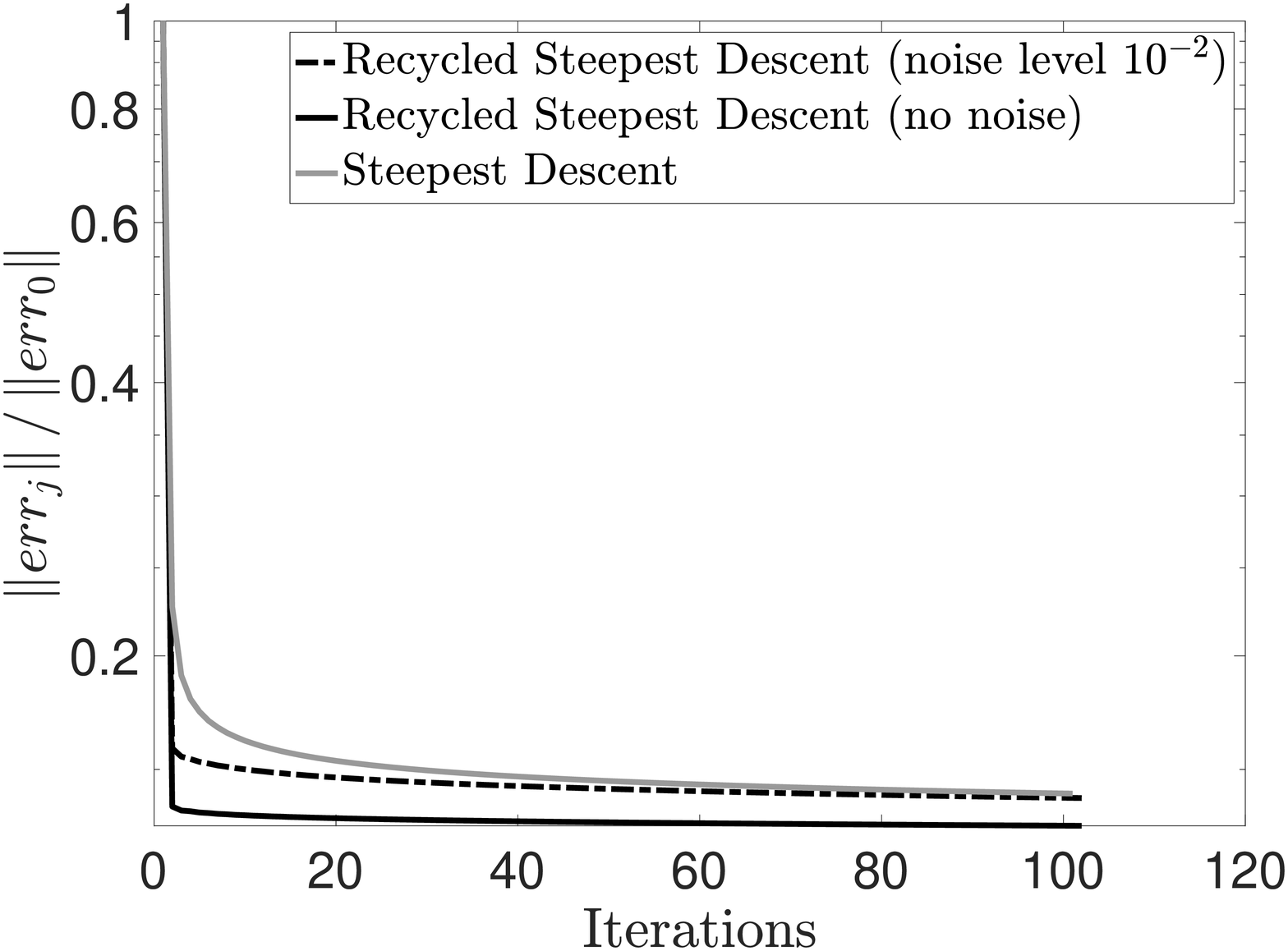}
	\caption{\color{black} Residual and error convergence for large-scale Gaussian blur problem with $\sigma=6$ and augmentation
	space generated using vectors resulting from applying 2 iterations steepest descent to the 
	noisy (dashed curve) and noise-free (solid black curve) data using
	Gaussian blur operators for five smaller standard deviations chosen equally-spaced from the interval
	$[0.5, 1.5]$.
	}
\end{figure}

{\color{black}
Our first experiment concerns an academic problem which serves
as a proof of concept. However, we still execute a large-scale version of this experiment using an
implicitly-defined, matrix free problem generated by the IR-Tools software package \cite{ir-tools}.  The IR-tools
package does provide an out-of-the-box generator function for the Gaussian blurring problem, {\tt PRBlurGauss()};
however, it does not allow one to specify precise standard deviation (i.e., the spread) of the point-spread function (PSF).
Looking under the hood, though, one sees that there is a helper function called {\tt psfGauss()} which generates  
a PSF function for a Gaussian with specified spread, which can then be fed into the general purpose
{\tt PRBlur()} problem generator, which outputs a matrix-free operator which convolves images with the PSF as well
as the true solution and noisy right-hand side.  A variety of solution images are available; we choose a predefined geometric
pattern for this experiment.  To make the problem sufficiently large-scale, we choose the image size to be 
$500\times 500$, meaning the implicit side of the matrix induced by the PSF is $250000\times 250000$. 
The IR-tools package is careful to not generate exact right-hand side data unless it is specified.  We rerun 
{\tt PRBlur()} with the option {\tt CommitCrime=`true'} to obtain $\vek b_{true}$, but this is only used to understand the 
effects of noise in generating the recycling vectors.
}  
In these experiments, we 
{\color{black} generated the problem}
for standard deviation 
$\sigma = 6$. 
	To construct an augmentation subspace $\CU$, we selected vectors resulting from applying 
	{\color{black} two} iterations of steepest descent to the right-hand side data but for 
	different Gaussian blurring operators with five
	smaller standard deviations chosen equally-spaced from the interval $[0.5, 1.5]$.  An orthonormal basis for the span of these (nine) vectors was obtained and taken as $\CU$.  {\color{black} We also generated a set of ``clean'' recycled vectors using the 
	noise-free $\vek b_{true}$, to see what effect noise has on the recycling process itself.
	It is demonstrated that such vectors do indeed offer better performance as recycled vectors, meaning that noise does
	effect the recycling process.}

	In \Cref{fig.augSteepestDescent-resid}, we see that for both problems, this strategy yields an improvement in convergence, both for the noise-free and noise-contaminated problems.  
	This is a promising numerical result and a proof of concept, but it requires
	further analysis {\color{black} and experimentation} to understand it in the proper concept to guide the use of this strategy for real-world, non-academic problems.

	\begin{figure}\label{fig.trueb-sdrecon}
		\centering
		\caption{The right-hand side and the best reconstruction produced by
		steepest descent for each, respectively.}
		\includegraphics[scale=.15,bb=0 0 1120 840]{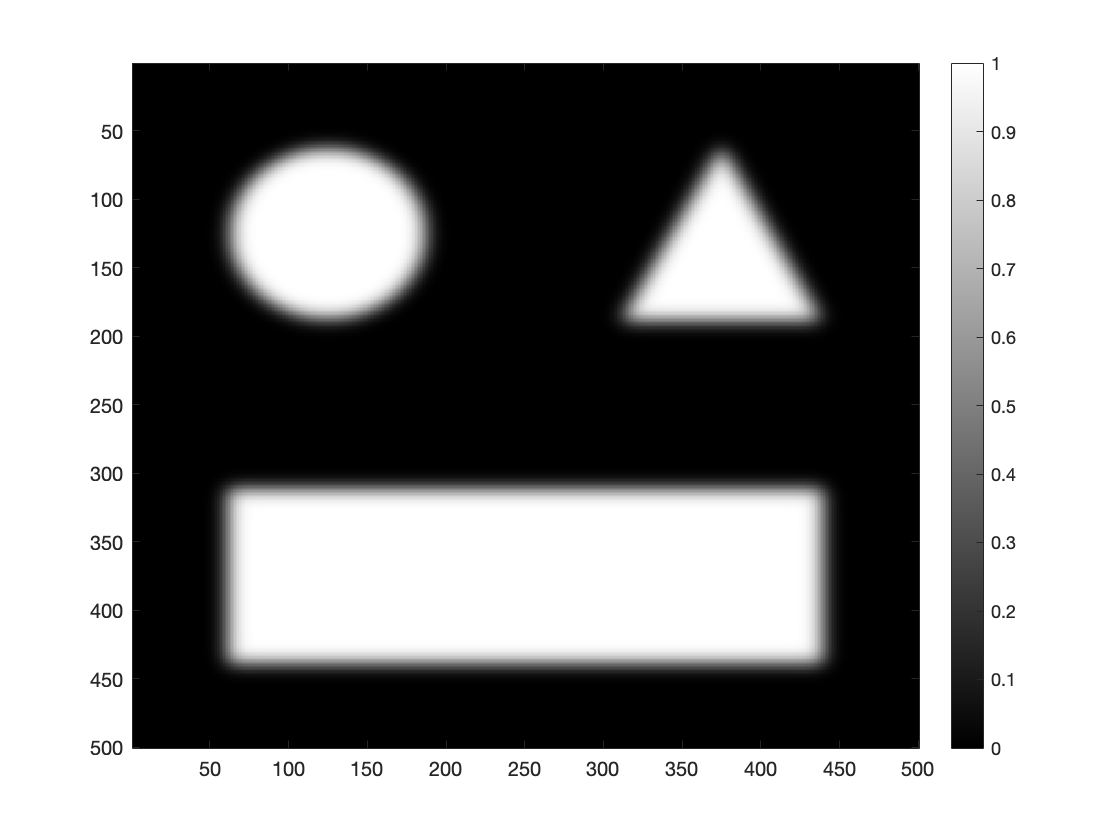}
		\includegraphics[scale=.15,bb=0 0 1120 840]{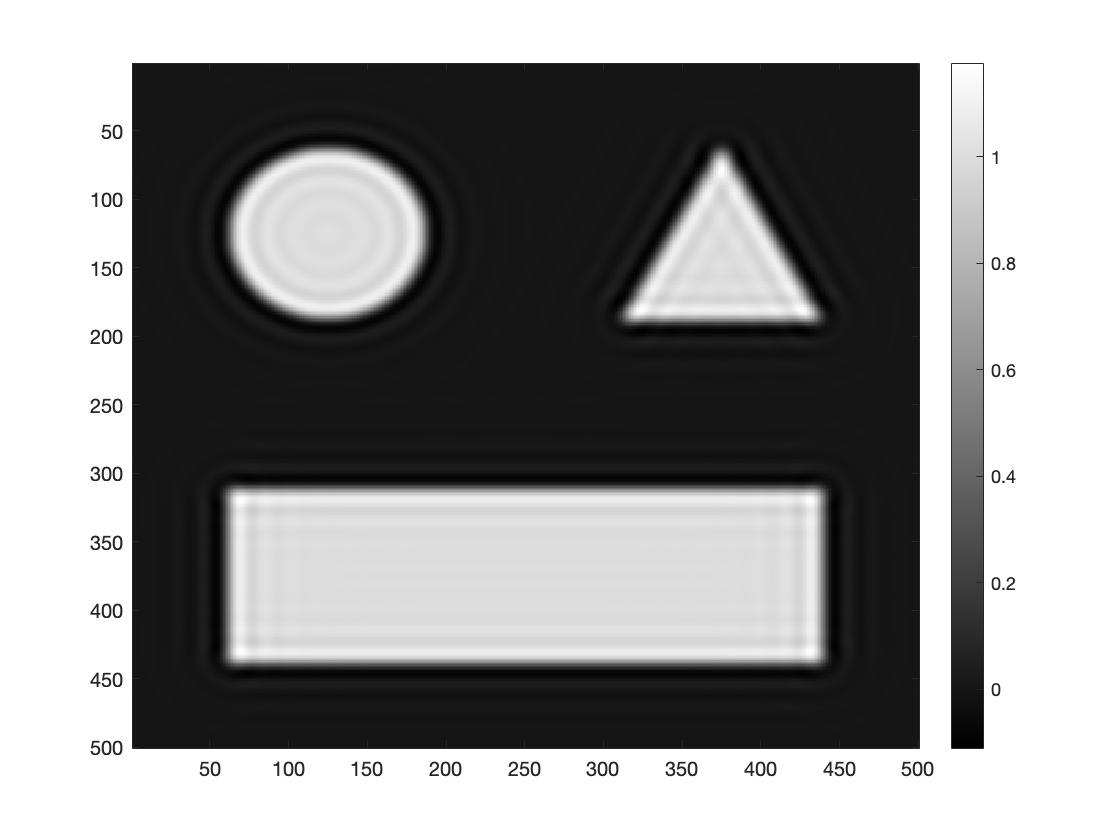}\\
		\includegraphics[scale=.15,bb=0 0 1120 840]{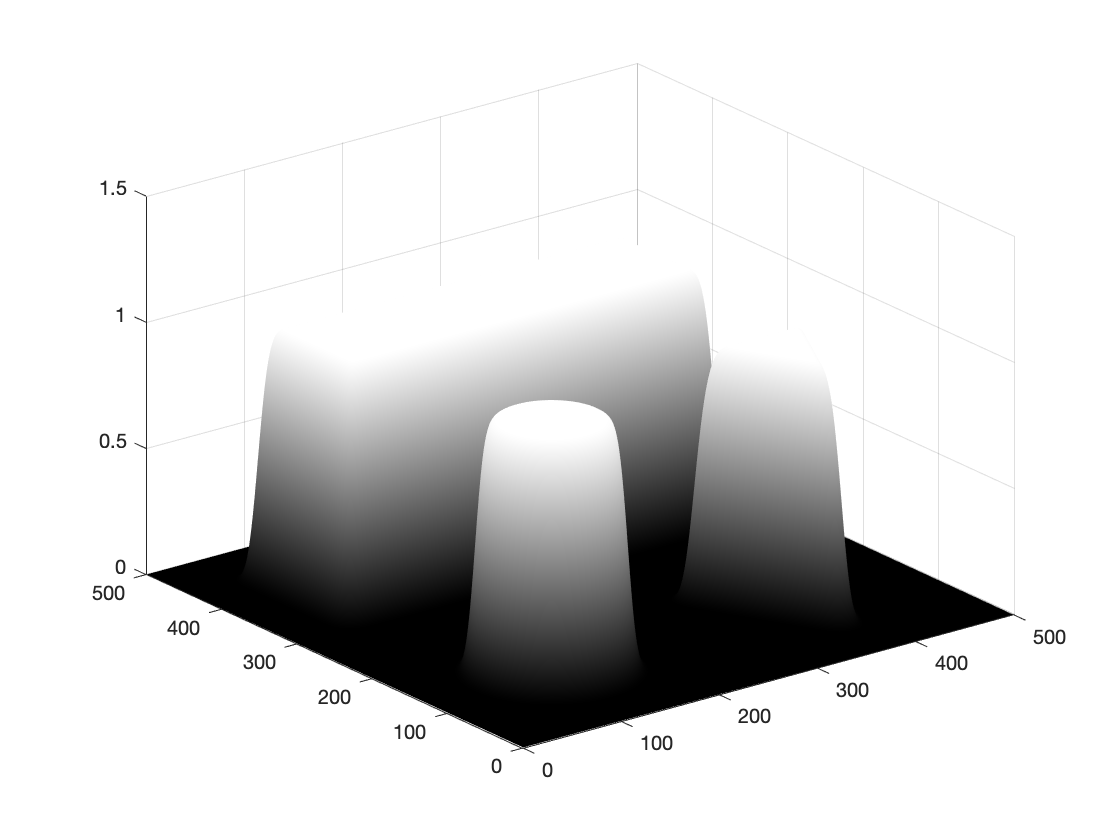}
		\includegraphics[scale=.15,bb=0 0 1120 840]{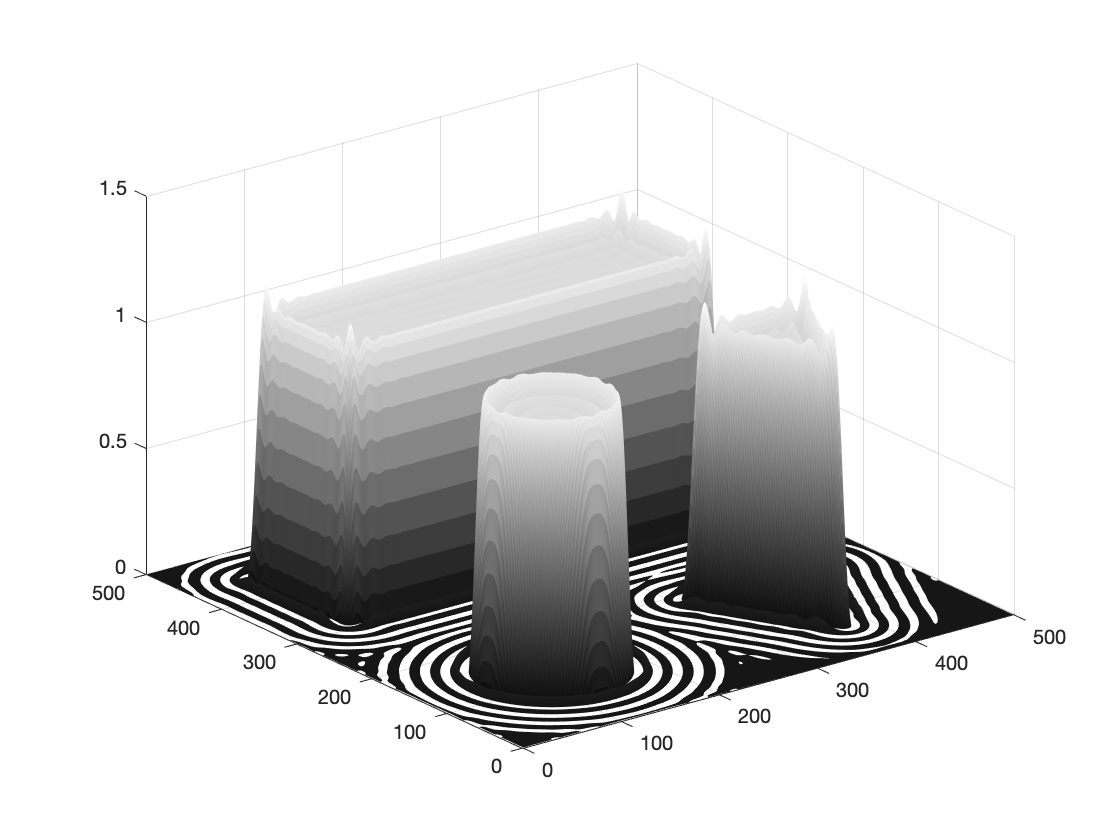}
	\end{figure}

	\begin{figure}\label{fig.augSd-recon}
		\centering
		\caption{The true image and reconstructed images using augmented steepest descent 
		 We show them both as images in the top row and
		three-dimensional surfaces in the bottom row.}
		\includegraphics[scale=.15,bb=0 0 1120 840]{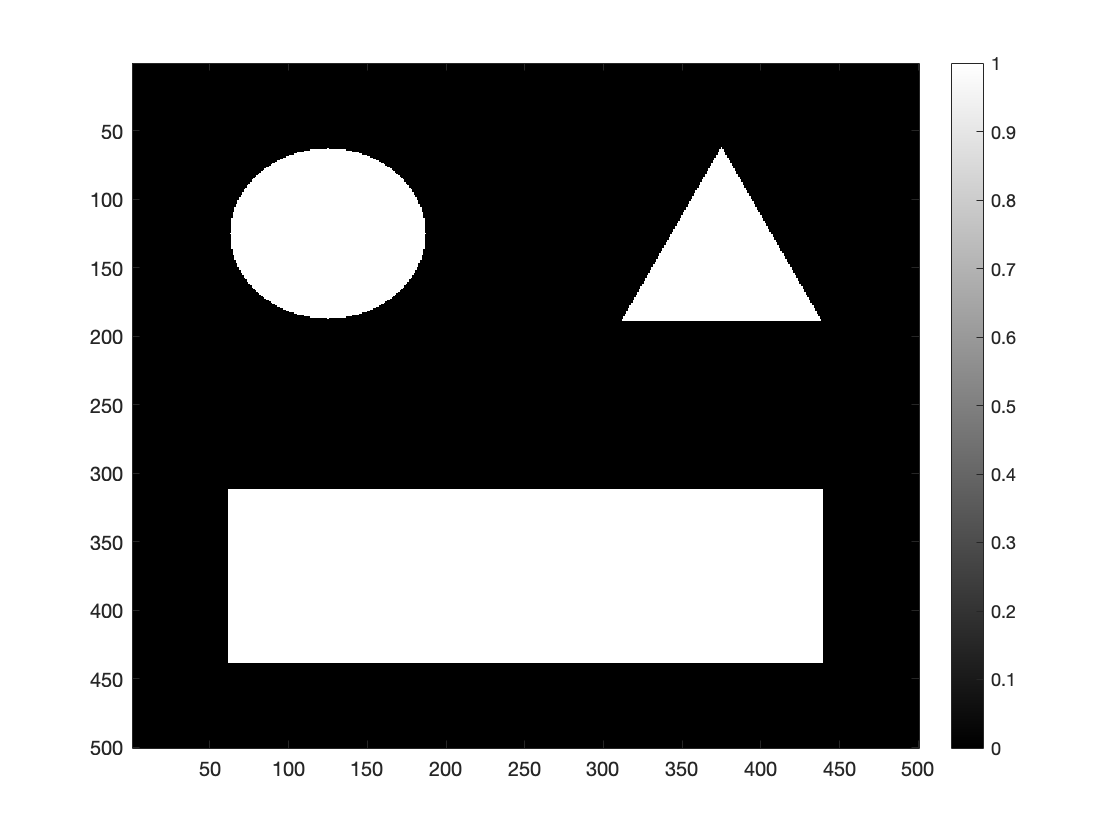}
		\includegraphics[scale=.15,bb=0 0 1120 840]{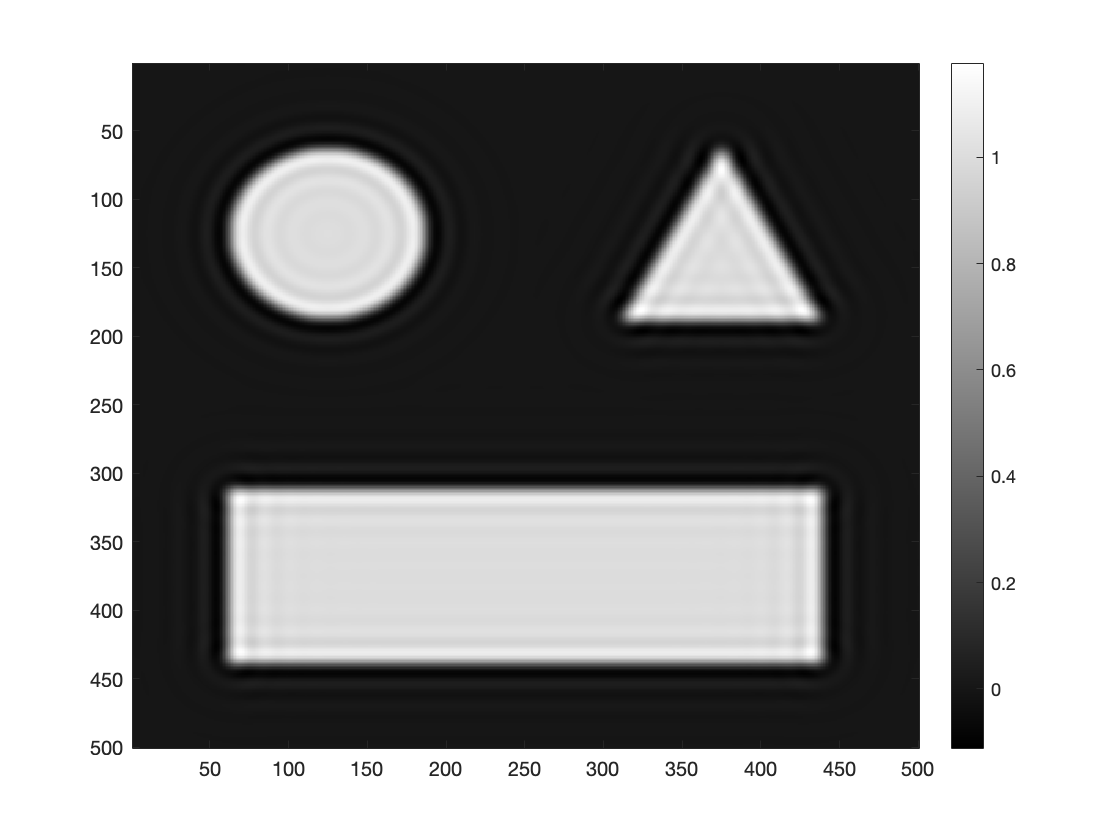}\\
		\includegraphics[scale=.15,bb=0 0 1120 840]{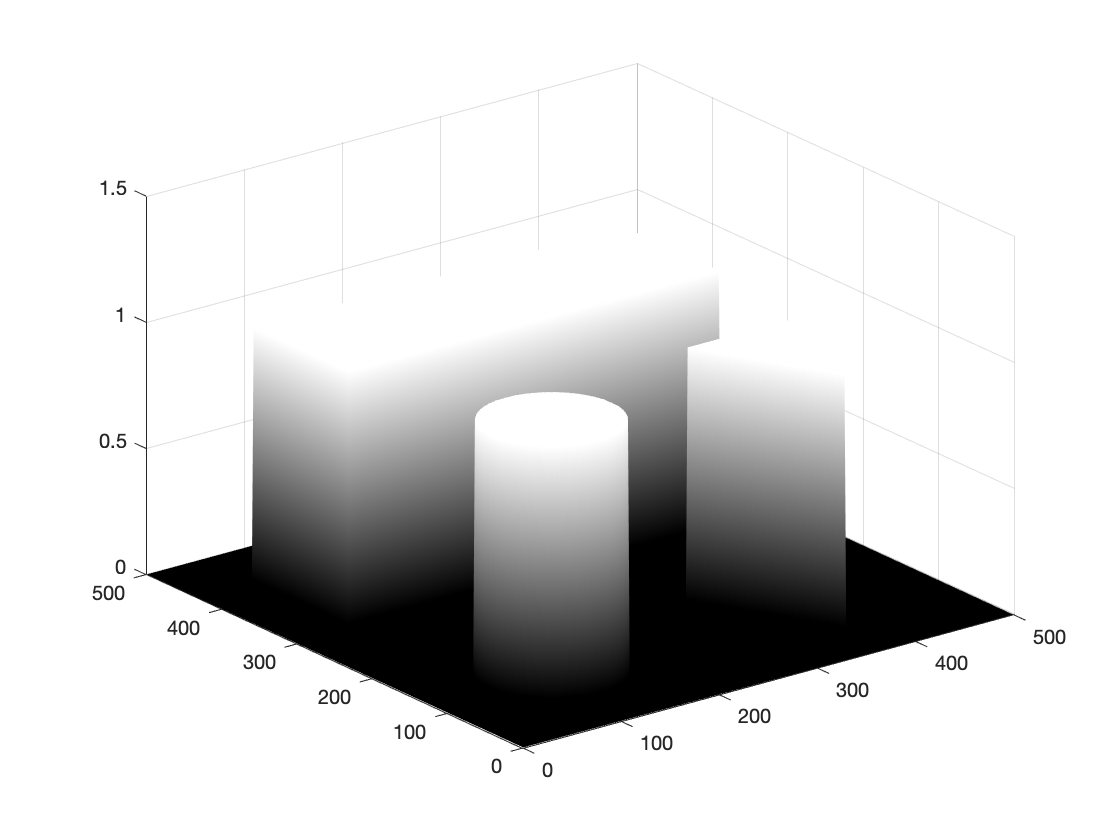}
		\includegraphics[scale=.15,bb=0 0 1120 840]{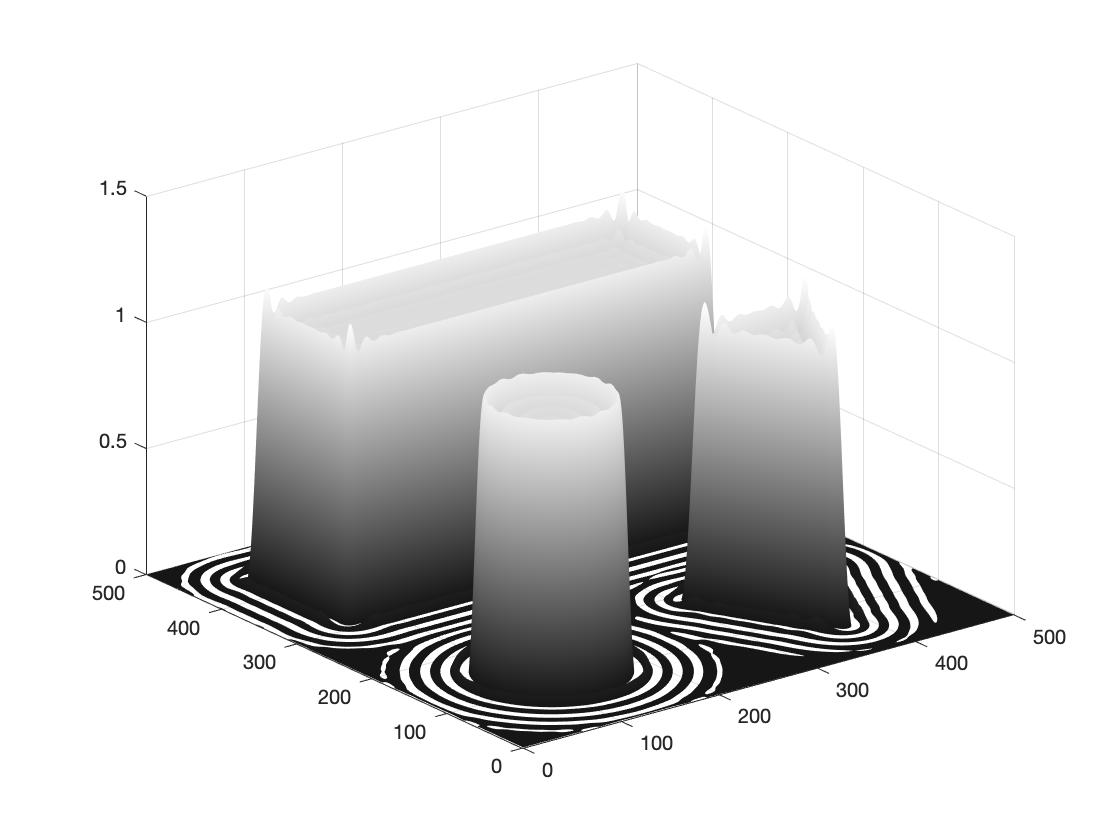}
	\end{figure}

\subsection{Adaptive optics: wavefront reconstruction}\label{AOsys}
\begin{figure}[H]
	\centering
	\includegraphics[scale=.35,bb=190 15 652 580]{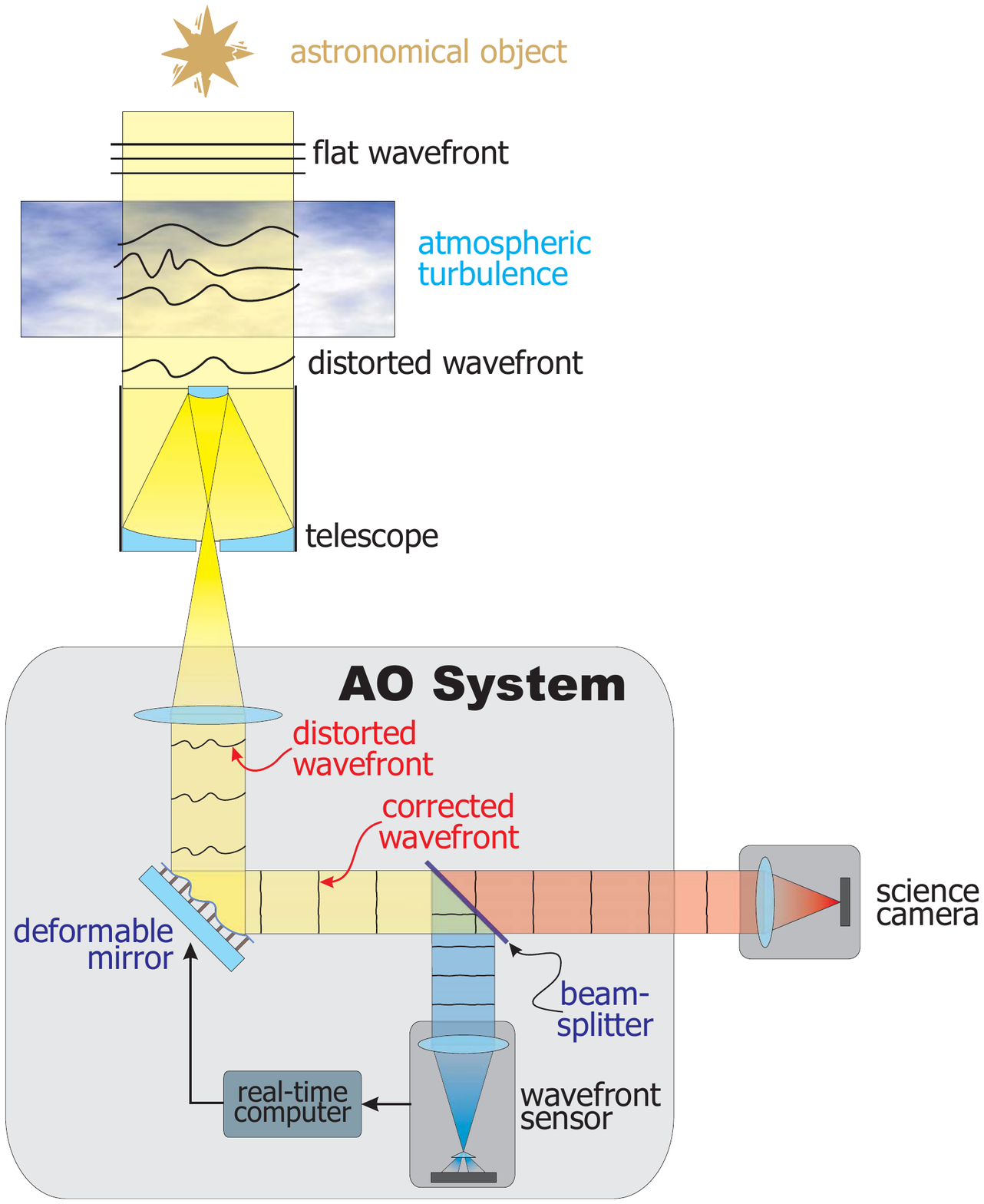}
	\caption{Principle of SCAO \label{fig.SCAO}}
\end{figure}

Turbulences in the atmosphere have a significant impact on the imaging quality of modern ground based telescopes. In order to correct for the impact of the atmosphere, {\it Adaptive Optics (AO)} systems are utilized. In these systems, the incoming light from bright guide stars is measured by wavefront sensors in order to obtain information of the actual turbulence in the atmosphere. In a next step, deformable mirrors use this information for an improvement of the obtained scientific images. For a review on the different mathematical aspects of AO systems we refer to \cite{ElVo09}. E.g., {\it Single Conjugate Adaptive Optics (SCAO)} systems use one sensor and one mirror for the correction and are able to obtain a good image quality close to the direction of the guide star. Figure \ref{fig.SCAO} shows a setup for a SCAO system. For the instruments of the new generation of telescopes, e.g., the Extremely Large Telescope (ELT) of the European Southern Observatory (ESO), the Pyramid sensor (P-WFS) will be used frequently to obtain the wavefront. A linear approximation of the connection between the incoming wavefront $\varphi$ and and the sensor measurements $(s_x,s_y)$ - assuming a non-modulated sensor - can be described as

\begin{equation}
	s_x (x,y) = \frac{1}{\pi}\int\limits_{-B(y)}^{B(y)} \frac{\Phi (x',y)}{x-x'}
\, dx' ~,
\end{equation}
analogous for $s_y$. An overview of existing reconstruction methods for the P-WFS can be found in \cite{ShaHuRam2020,Hut17,Hut18}, and a throughout analysis of the various P-WFS models is given in \cite{HuShaRa19_1,HuShaRa19_1}.

\begin{figure}[H]
	\centering
	\includegraphics[scale=.5,bb= 0 0 679 459]{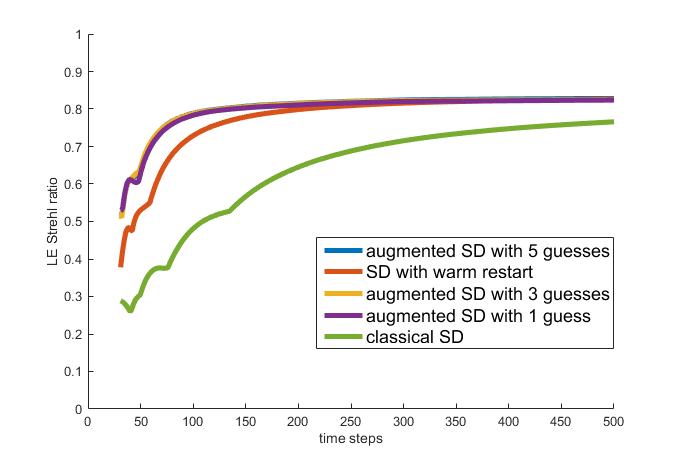}
	\caption{Monitoring of wavefront reconstruction accuracy over time. \label{fig.vici}}
\end{figure}

The linear models can be used in case of small wavefront aberrations which is the case in a {\it closed loop} setting, i.e., if the wavefront sensor is optically located {\it behind} the deformable mirror.
Due to the rapidly changing atmosphere, the shape of the deformable mirror has to be readjusted every 1-2 milliseconds over the whole imaging process, which might last several minutes. Therefore, the reconstructions of the wavefronts have to be accurate and have to be obtained in real-time.
As a quality measure we use the Strehl ratio ($SR$), which relates the imaging quality of the telescope without disturbing atmosphere ($SR = 1$) to the corrected imaging quality, i.e.. $SR\in [0,1]$. In the beginning of the observation the Strehl ratio is low, as the deformable mirror is not yet adjusted to the turbulence conditions. After a few time steps, the mirror is in a state where only small adjustments are necessary - {\it the loop has been closed}. This happens in our case when $SR\in [0.6, 0.7]$. In this state, the linear modes for the P-WFS are valid.
We show the utility of the recycled steepest descent method for multiple steps of an closed-loop iteration for a wavefront reconstruction for the instrument METIS of the ELT. The recycling subspaces were formed by the previous $1-3$ reconstructions, see \cref{fig.vici} for the results. In the long run, both steepest descent with recycling and warm restart (i.e, where the last reconstruction is used as a starting value for the new reconstruction) yield the same Strehl ratio. However, the method with recycling is able to close the loop much faster, which is important for a stable and efficient imaging of the scientific object.

\subsection{Adaptive optics: Image reconstruction}
Our final example concerns another problem related to astronomical imaging. As for any optical system, the design of the system determines the obtainable imaging quality: Instead of the "true" image one always obtains a blurred image. The extend of blurring is described mathematically by the point spread function (PSF) of the optical system, and the measured image can be modeled as the convolution of the true image and the PSF:

\begin{equation}\label{blurring}
	I(t) = \int I_{true}(s)\cdot PSF(t-s)\, ds.
\end{equation}
Thus, the original image can be recovered - assuming the PSF is known - by solving the linear ill-posed problem \eqref{blurring}.
In this example we want to reconstruct an astronomical image from the measured image obtained by a telescope
using the approximate  PSF of the telescope. It is well known that the PSF of a telescope without a turbulent atmosphere above the telescope (e.g., a telescope in space) is given by
\be\label{eqn.PSF}
	PSF_{\color{black}0}(x) = \ab{\CF \{P\}(x)}^{2}
\ee
where $P$ is the characteristic function of the aperture of the main mirror of the telescope (usually an annulus) and $\CF(\cdot)$ represents the Fourier transform of the argument. {\color{black}Please note that the recovery of the original image from the measured one results in most cases in an ill posed problem, at least for smooth PSF, which is, e.g., the case for \eqref{eqn.PSF}.} For ground based telescopes, however, the turbulent atmosphere {\color{black} above the telescope }has to be taken into account. {\color{black} The related PSF can } be modeled as
\be\label{eqn.PSF-def}
	PSF_{\phi}(x,t) = \ab{\CF \{P(\cdot)e^{i\phi(\cdot,t)}\}(x)}^{2}.
\ee
Here, $\phi(r ,t)$ describes the turbulence above the telescope. {\color{black}As we have seen above, an AO systems aims to reduce the impact of the turbulence to the imaging process - see Section \ref{AOsys}. For various reasons, e.g., the time delay between measuring the turbulences and their correction or the limited resolution of wavefront sensors and deformable mirrors, such an correction will never be perfect and therefore there will always remain a residual (uncorrected) turbulence. Based on the data from the AO system it is, however, possible to estimate those uncorrected turbulences and obtain an approximation of $PSF_{\phi}$ of the observation.}
A full description of this problem, {\color{black}i.e. the reconstruction of the PSF from AO data of the scientific observation}, can be found, e.g., in \cite[Section 5.4]{ReRaSoWa} and references therein.
As astronomical images are acquired over a period of time, one actually {\color{black}recovers a} time-averaged PSF, denoted $\ip{PSF_{\phi}(\cdot,t)}$. \\
{\color{black} In our experiments, the time averaged $PSF_\phi$ with residual turbulences taken into account is simulated by OCTOPUS \cite{LVK06,LeLouarn_OCTOPUS_04}, the simulation tool of the European Southern Observatory.}
In addition, we use a representation of $PSF(x)$, the PSF of the telescope without turbulences. Both PSFs are represented by matrices, i.e., as the functions
spatially discretized for images of a certain size. {\color{black} Please note that the recovery of the original image from the convolution data remains ill-posed in both cases.}

\begin{figure}
\centering
\caption{The true and blurred {\tt Star Cluster} image \label{fig.star-cluster} with generated stars of various brightnesses,
with the true image shown in log-plot to clearly show
all the generated stars clearly}
		\includegraphics[scale=.15,bb=0 0 1120 840]{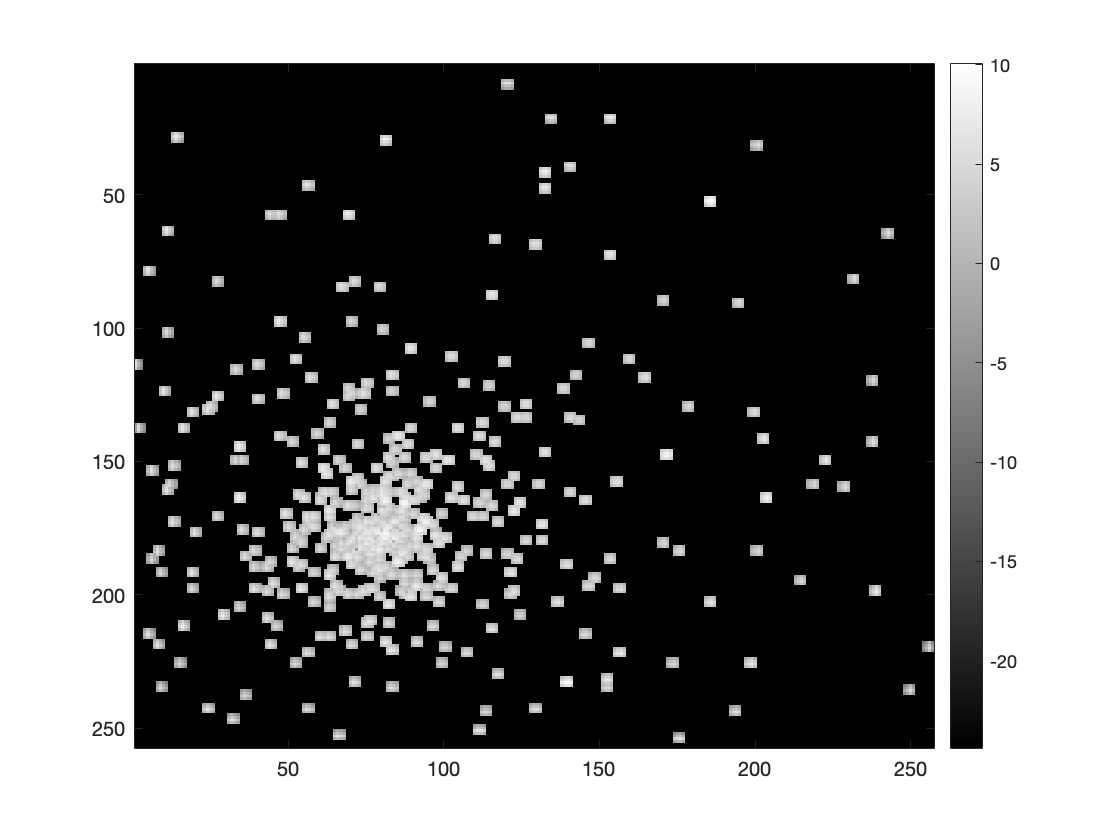}
		\includegraphics[scale=.15,bb=0 0 1120 840]{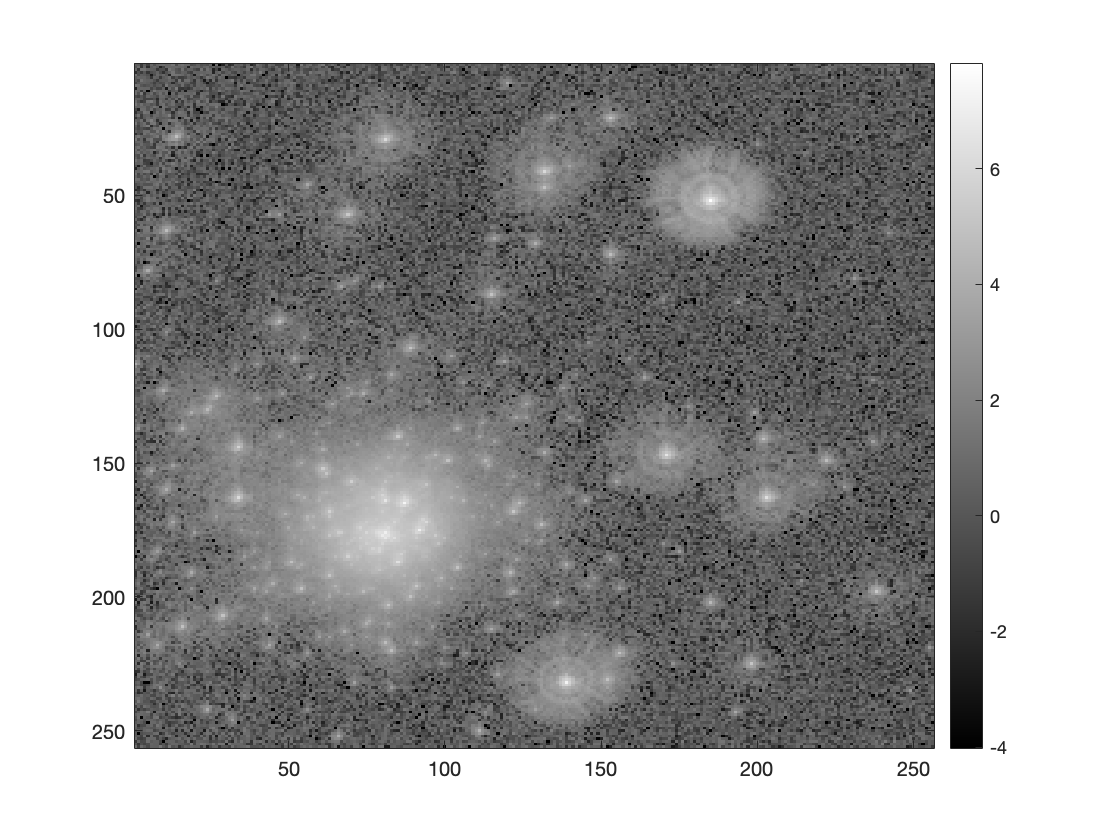}
\end{figure}

Convolution of images with PSFs was performed using the function {\tt psfMatrix()} available, .e.g., as a part of the {\tt Matlab}
toolbox {\tt IRTools} \cite{ir-tools}.  We use the {\tt Star Cluster} image, which can be obtained from the test problems of the
{\tt Matlab} package {\tt Restore Tools} \cite{restore-tools}.  The provided artificial {\tt Star Cluster} image is $257\times 257$; so for these
experiments, the discretized PSFs are resized (i.e., downscaled) to also be $257\times 257$ using the {\tt imresize()} function of {\tt Matlab}.
A {\tt psfMatrix} object is then instantiated which can be applied via multiplication to an image to perform a convolution, i.e.,
\be\nn
	\tt Apsf = psfMatrix(psf\_257,`periodic',[C\ C],[N\ N])
\ee
where $N=257$ and $C=129$, the center coordinate of the image, and we use periodic boundary conditions.  Entry $(C,C)$ of the PSF indicates
where the value of the PSF at $(0,0)$ is located, where we assume the square image is centered at the origin.  The PSF inducing the adjoint
operator can be created via the relation $PSFadj(x,y) = PSF(-x,-y)$ which, along with the image being centered at the origin, can be used to
construct the discrete PSF inducing the discretized adjoint operator \cite{roland-email}.
In \Cref{fig.star-cluster}, we show log plots of the true {\tt Star Cluster} image and its blurred counterpart convolved with the PSF representing the
atmospheric distortion.  The blurred image was perturbed with random noise generated with uniform distribution using {\tt Matlab}'s {\tt rand()}
function, which was then scaled to have norm $10^{-1} \nm{b}_{F}$, where $b$ represents the convolved true image shown on the right in
\Cref{fig.star-cluster}.

{\color{black} In this study, we demonstrate that subspace recycling works well for the considered deconvolution problem. As for a telescope with a good AO system $PSF_0$ is at least a coarse approximation of $PSF_{\phi}$, we propose to form the recycling subspace from a few eigenvectors of the matrix representing $PSF_0$.
Here, we take advantage of the fact that although $PSF_\phi$ changes for every observation,  $PSF_0$ remains constant.
This enables us to pre-calculate eigenvectors of its induced convolution operator ahead of time offline.}

\begin{figure}
\centering
\caption{Residual and error convergence of steepest descent (shown in red) and augmented steepest descent (shown in shades of gray)
with augmentation spaces being eigenvectors
of the operator induced by the telescope PSF, where we used the eigenvectors associated to the largest eigenvalues.
The curve's darkness is determined by how many eigenvectors were used in the
augmentation.  More vectors produces a darker curve. \label{fig.AOPSF-augSD-eig}}
		\includegraphics[scale=.15]{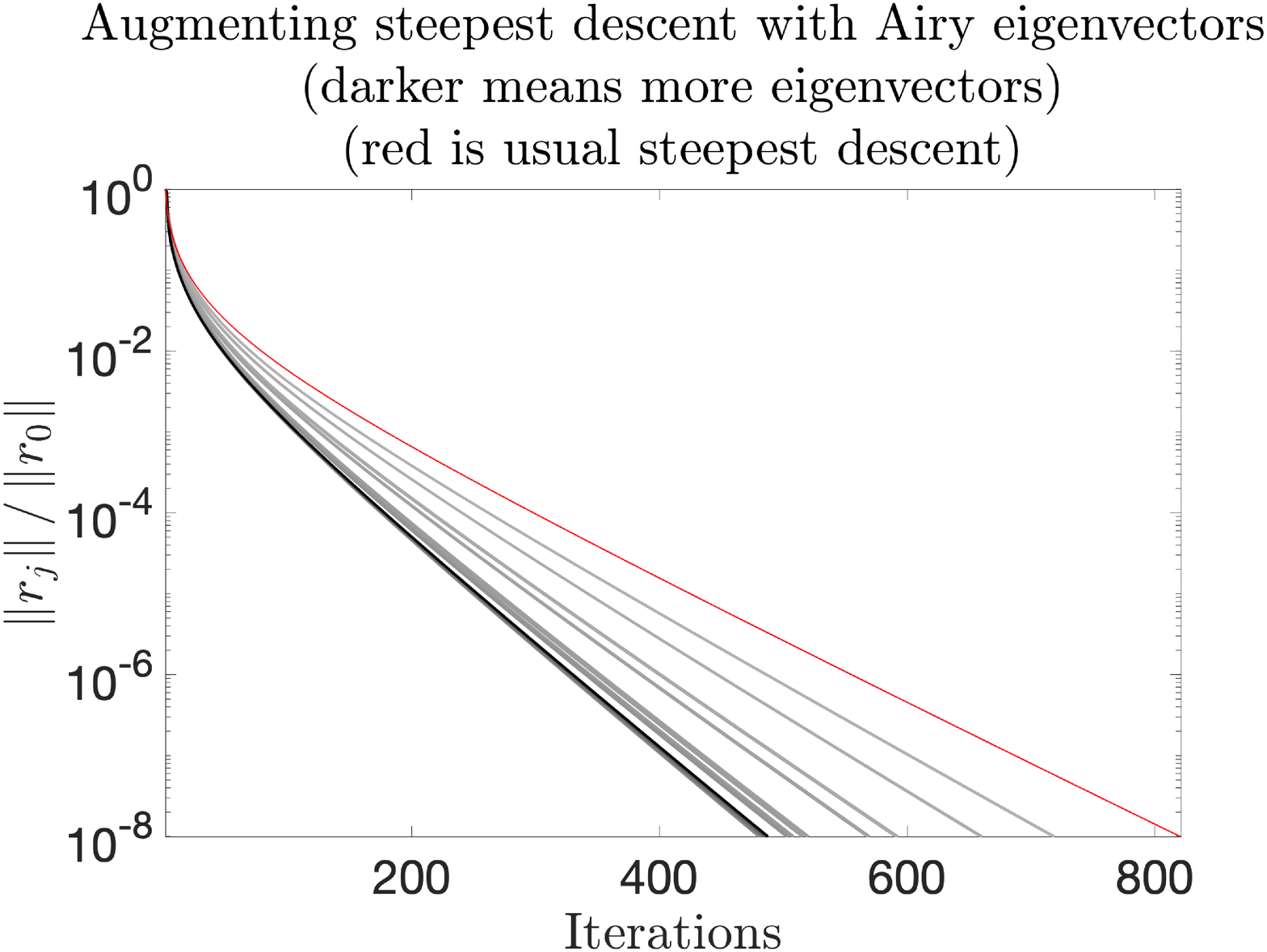}
		\includegraphics[scale=.15]{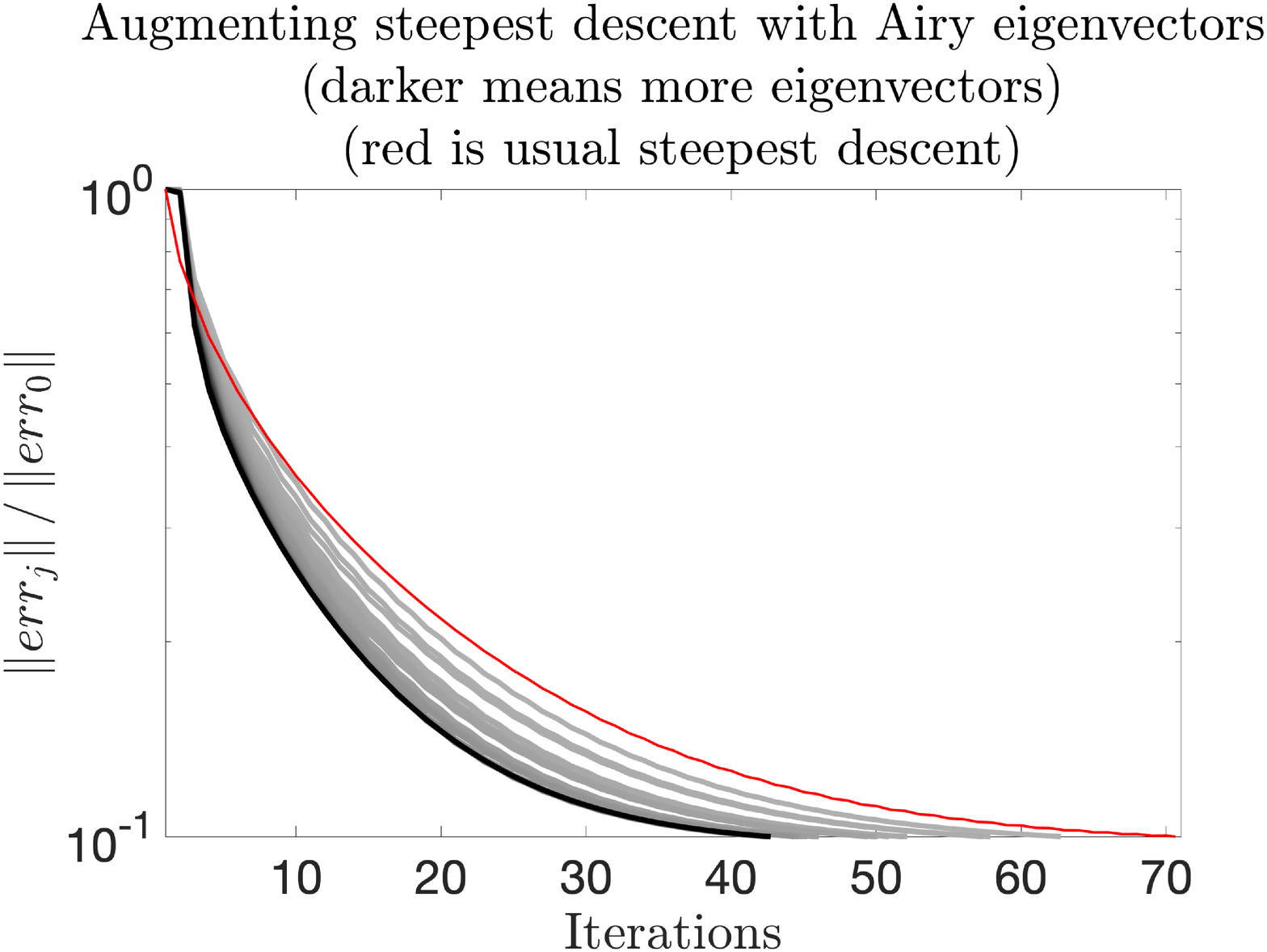}
\end{figure}

{\color{black}We calculated the telescope} eigenvectors via the matrix-free {\tt eigs()} routine of {\tt Matlab} which uses a Krylov-based iteration to compute eigenvalues and eigenvectors
for a few of the largest eigenvalues of the operator.  In this experiment, we chose to compute the top 200 but had to perform post-processing step to eliminate the spurious complex eigenvalues/-vectors that are created in pairs, leaving us with
198 {\color{black}eigenvalues}.
Ordering the eigenvalues thereof in
descending order, we augment with different collections of these pre-computed eigenvectors, starting with just the first, then the first two, and so on,
up to all 198. In \Cref{fig.AOPSF-augSD-eig}, we show the residual and error curves for all experiments together.  One sees that augmenting
with these eigenvectors is effective in increasing the speed of convergence of the steepest descent method, with acceleration increasing as
we add more vectors, but with diminishing returns.

\begin{figure}
	\centering
	\caption{We show (left figure) for each recycled subspace dimension from this experiment at
	what iteration semiconvergence and the discrepancy principle
	are achieved, and we also mark where jumps in the eigenvalues of the operator induced by the PSF (the Airy function) of the telescope occur.
	On the right, we show the eigenvalues of this operator.  \label{fig.conv-tracker}}
		\includegraphics[scale=.15]{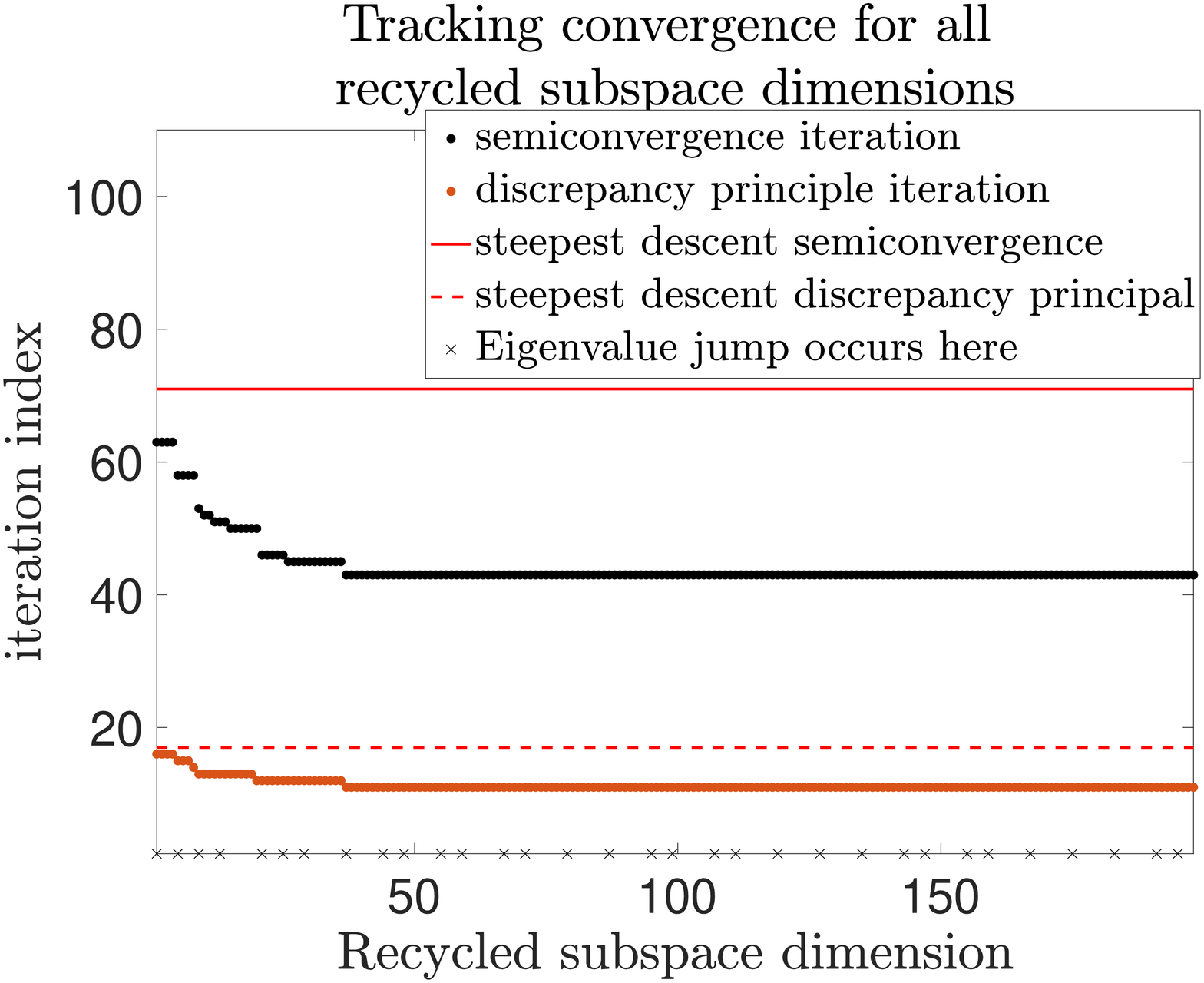}
		\includegraphics[scale=.15]{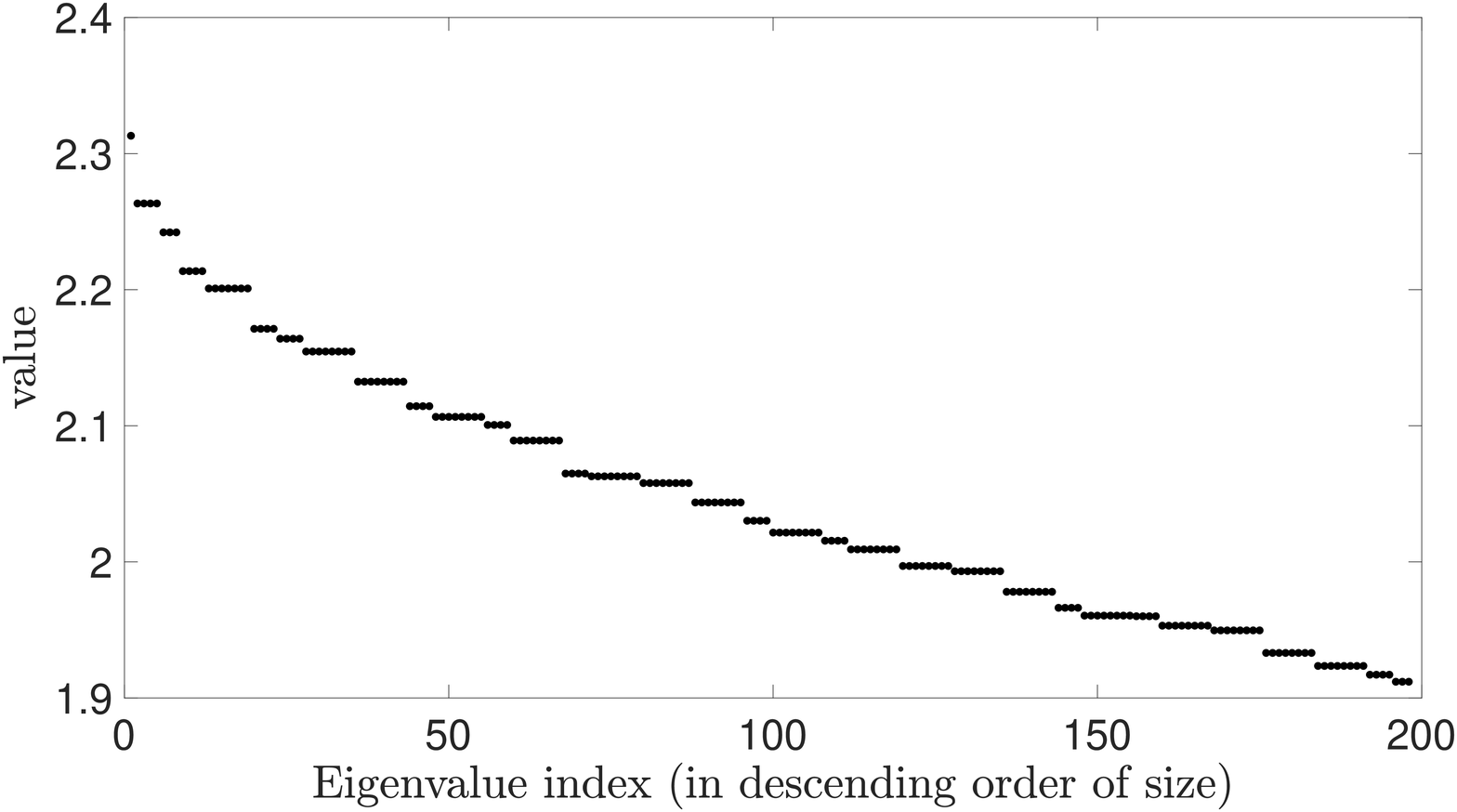}
\end{figure}

\begin{figure}
\centering
\caption{The first nine eigenfuntions (with eigenvalues in descending order) of the operator induced by the PSF, i.e.,
\cref{eqn.PSF-def}, of the telescope.  Note
that the first function is actually a constant up to roundoff noise. \label{fig.eigenvectors-airy}}
		\includegraphics[scale=.06,bb=0 0 2000 1480]{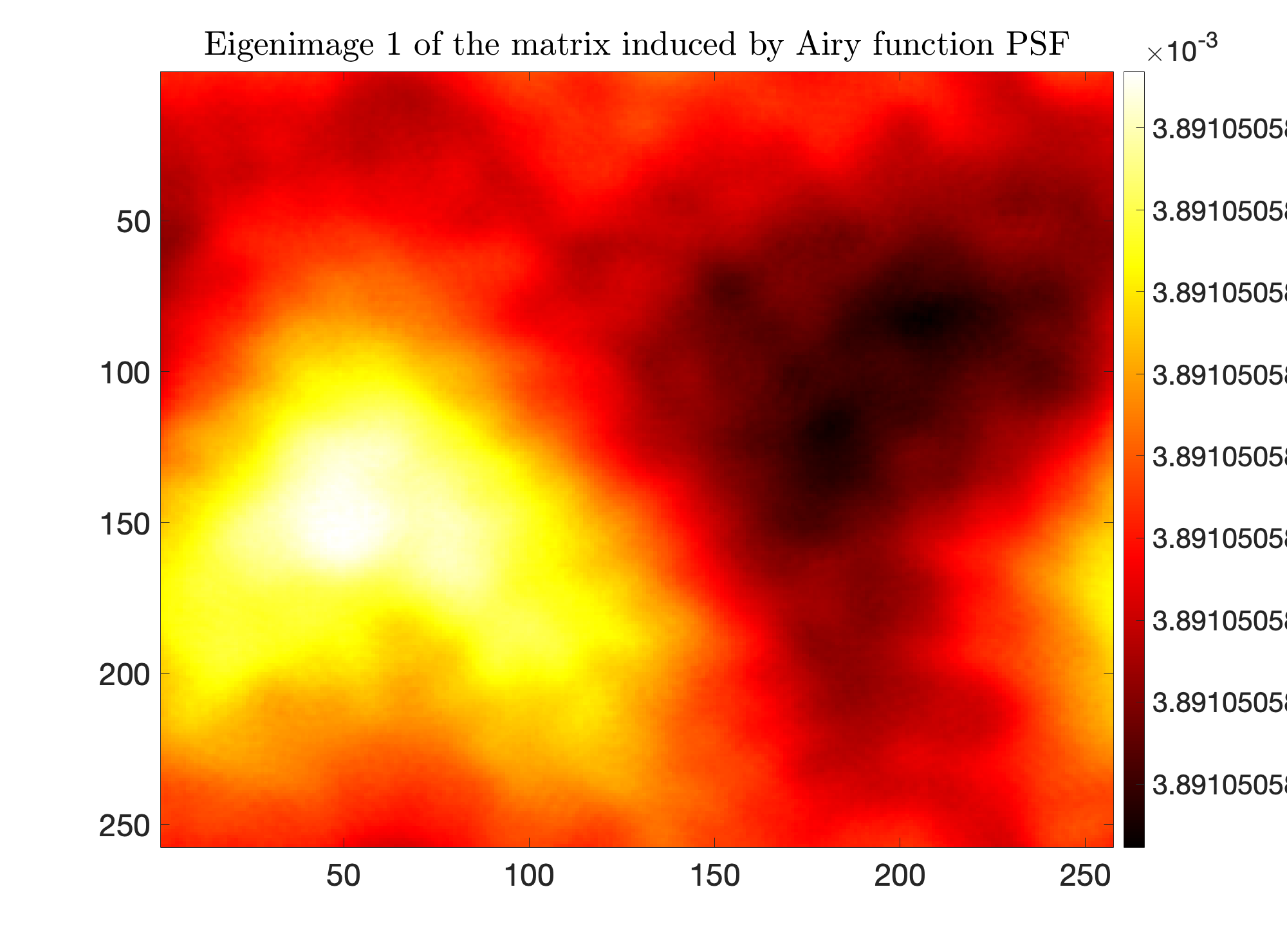}
		\includegraphics[scale=.06,bb=0 0 2000 1480]{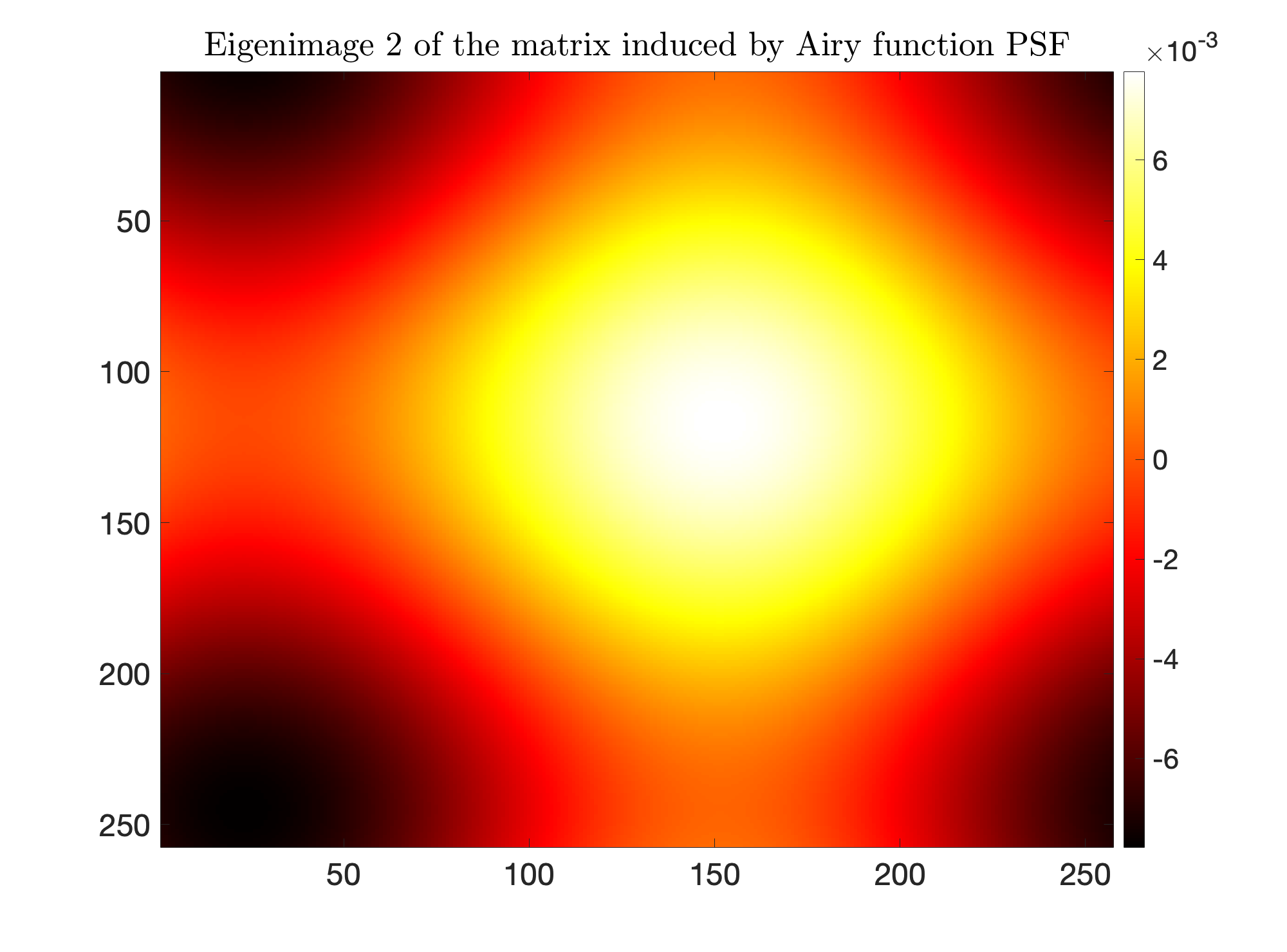}
		\includegraphics[scale=.06,bb=0 0 2000 1480]{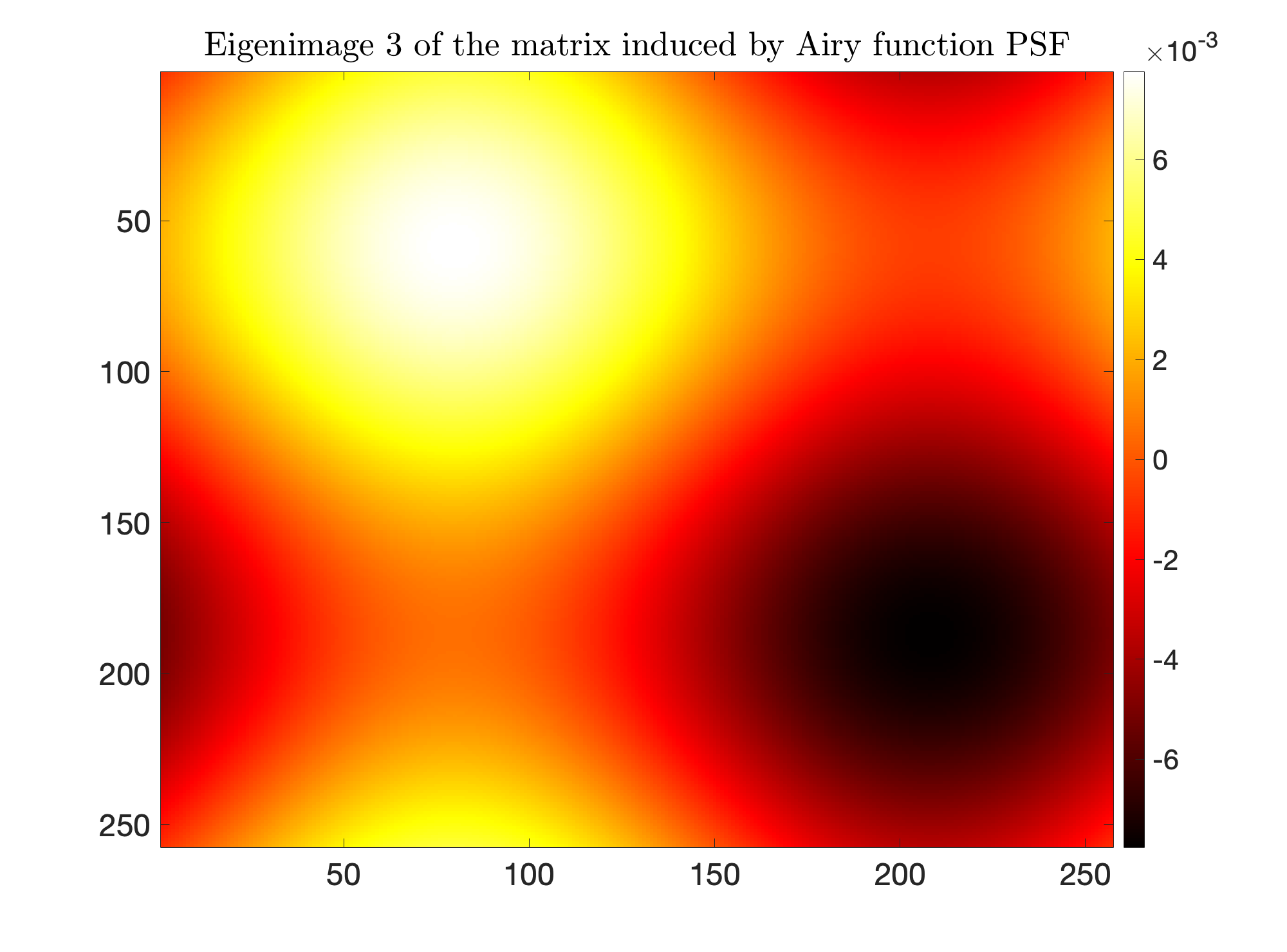}\\
		\includegraphics[scale=.06,bb=0 0 2000 1480]{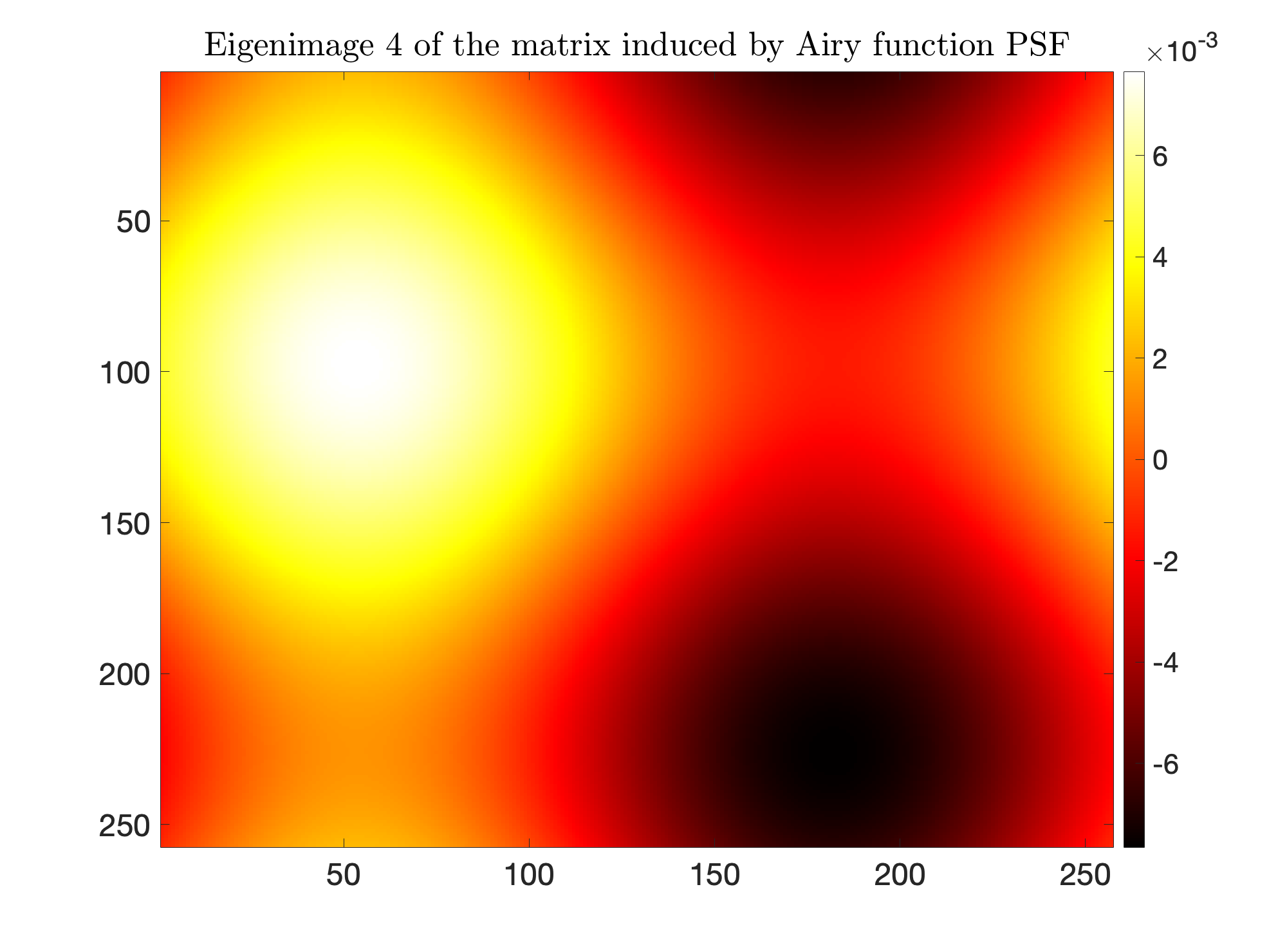}
		\includegraphics[scale=.06,bb=0 0 2000 1480]{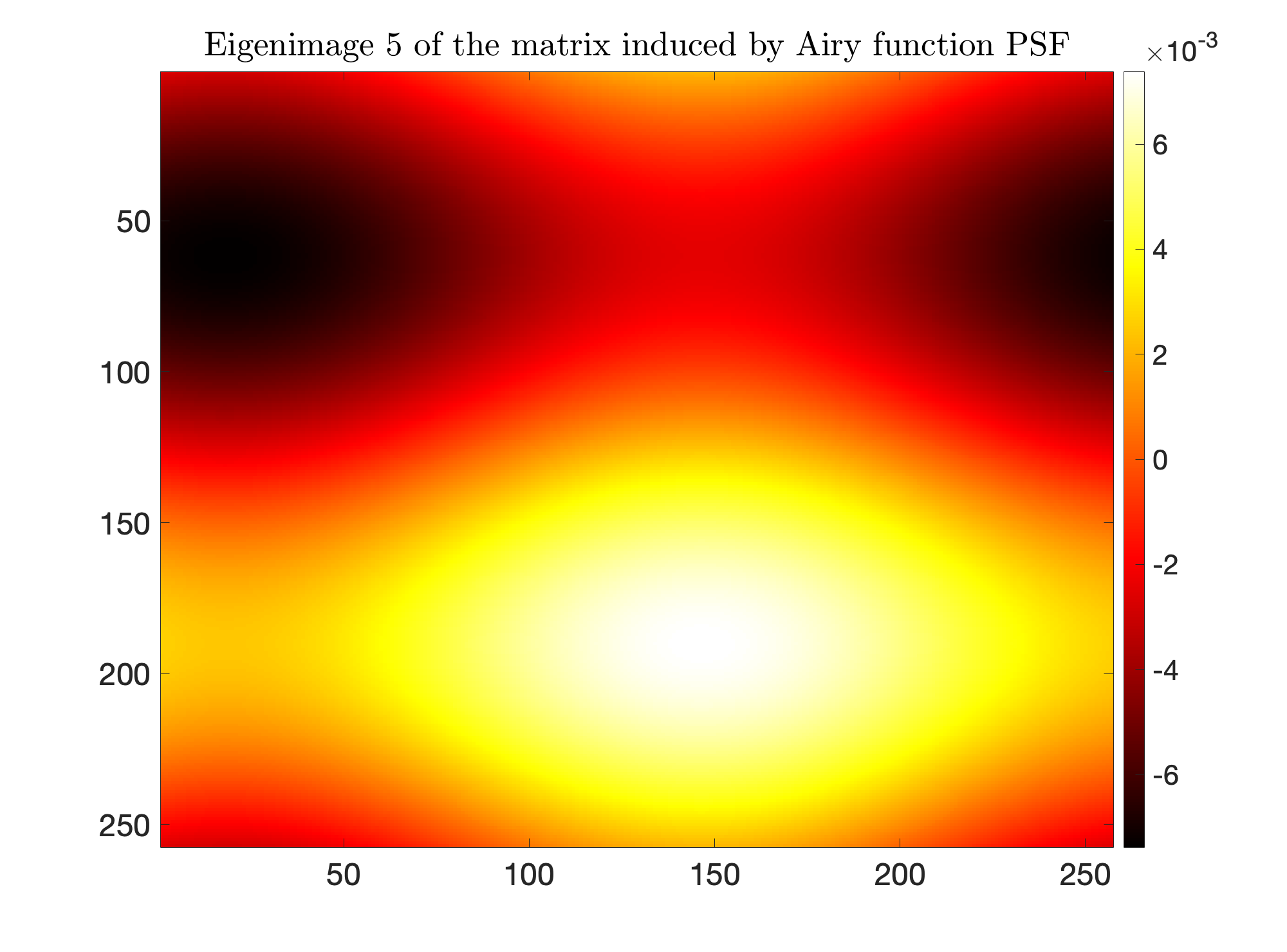}
		\includegraphics[scale=.06,bb=0 0 2000 1480]{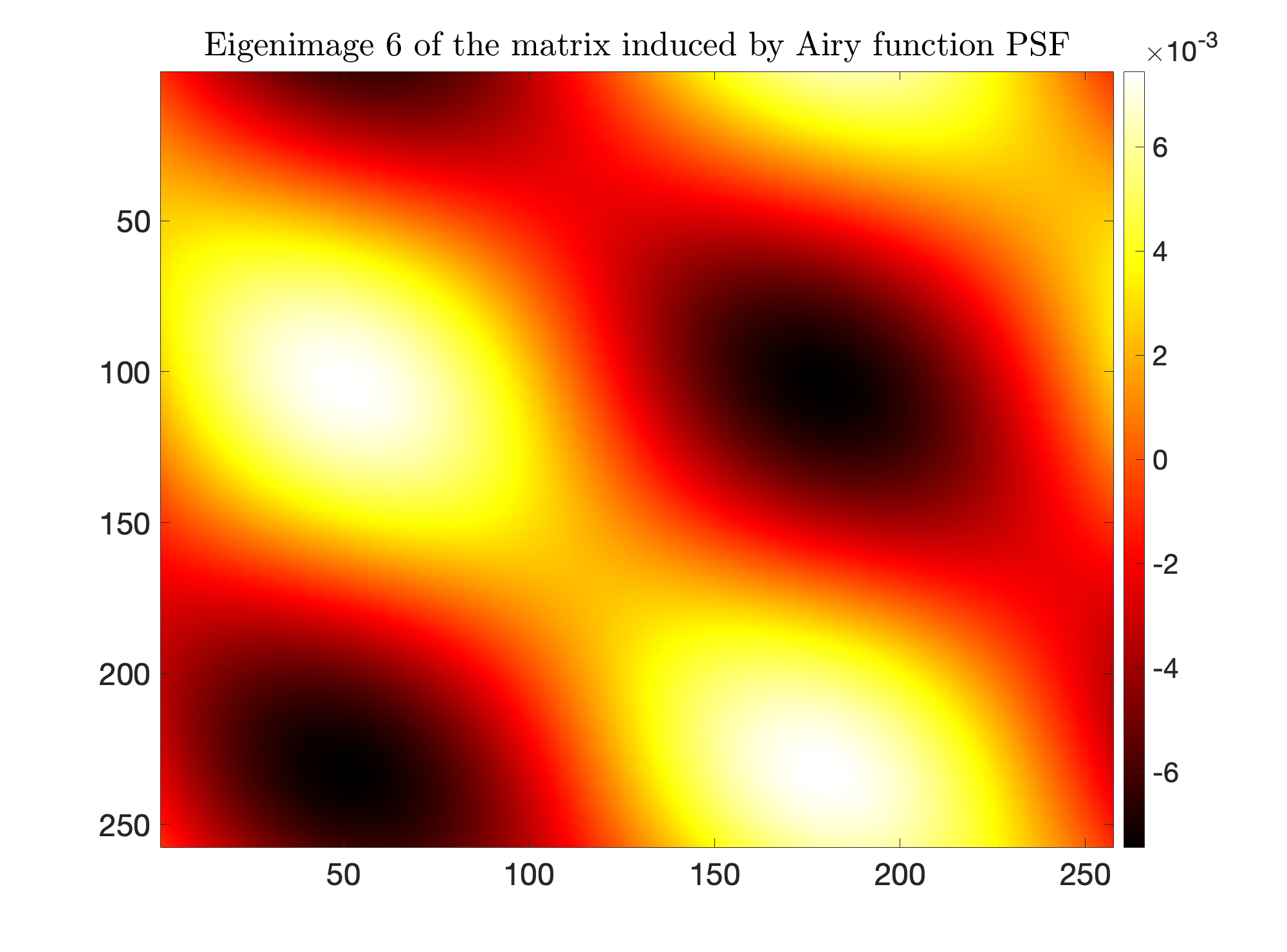}\\
		\includegraphics[scale=.06,bb=0 0 2000 1480]{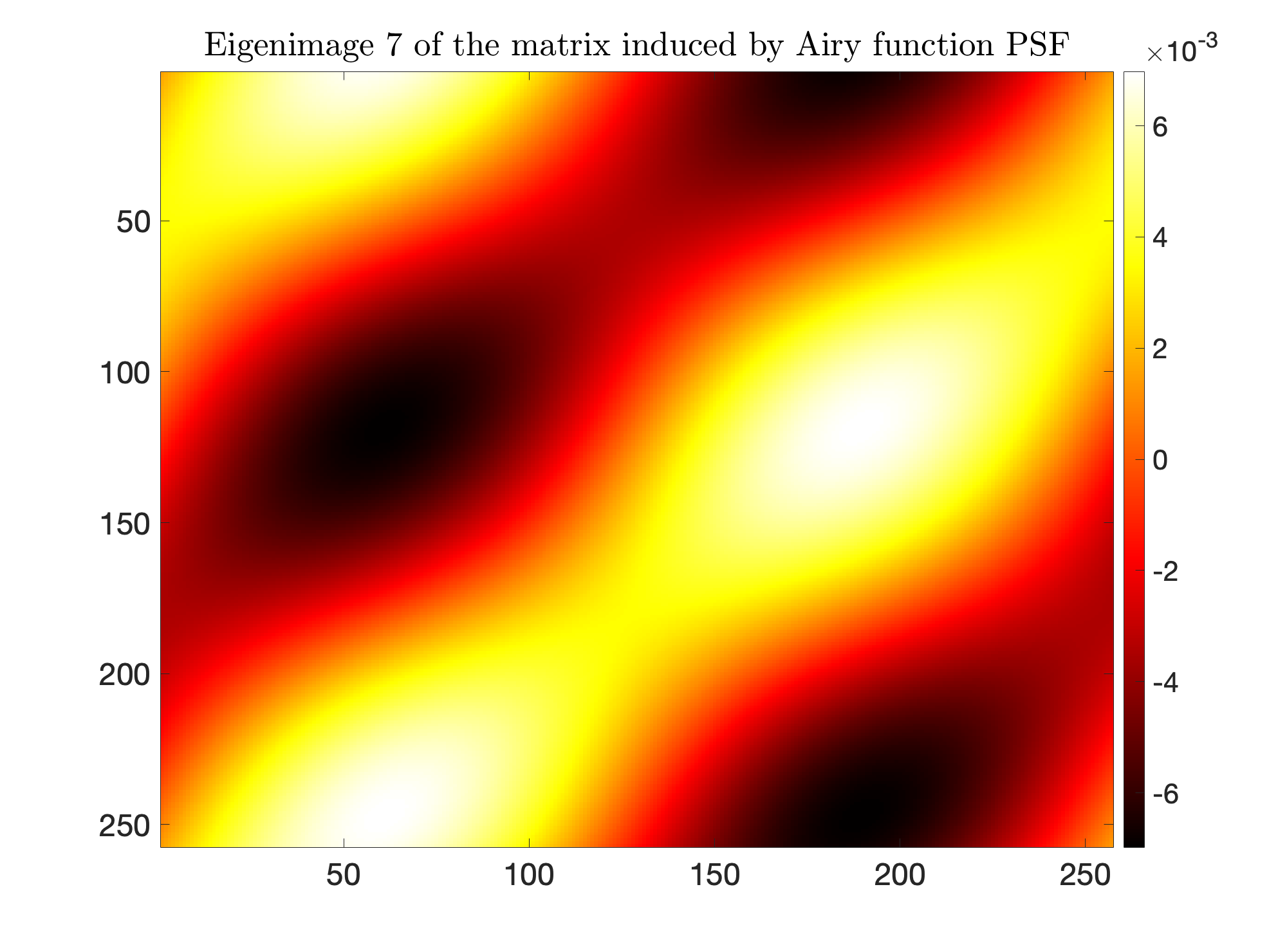}
		\includegraphics[scale=.06,bb=0 0 2000 1480]{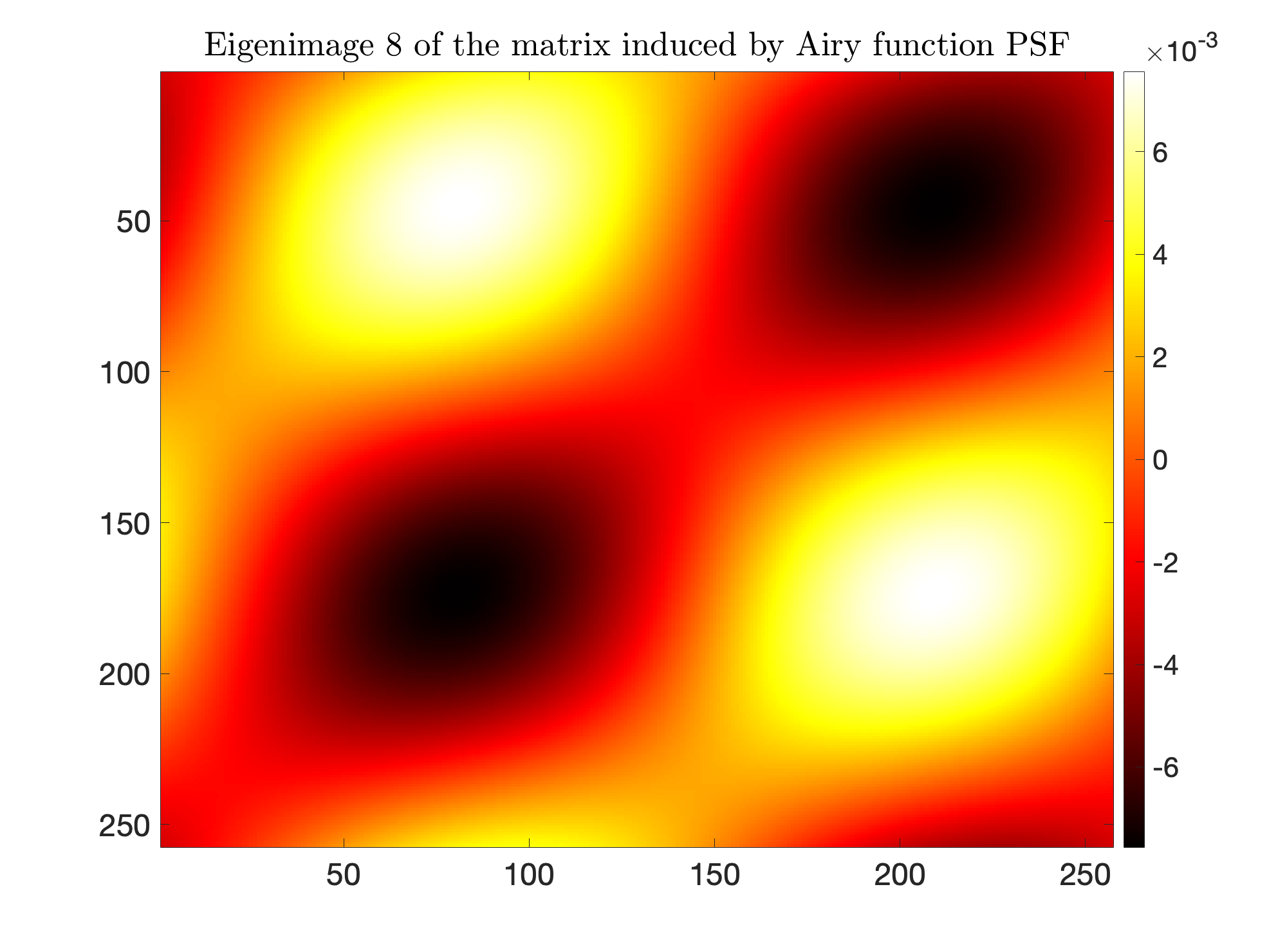}
		\includegraphics[scale=.06,bb=0 0 2000 1480]{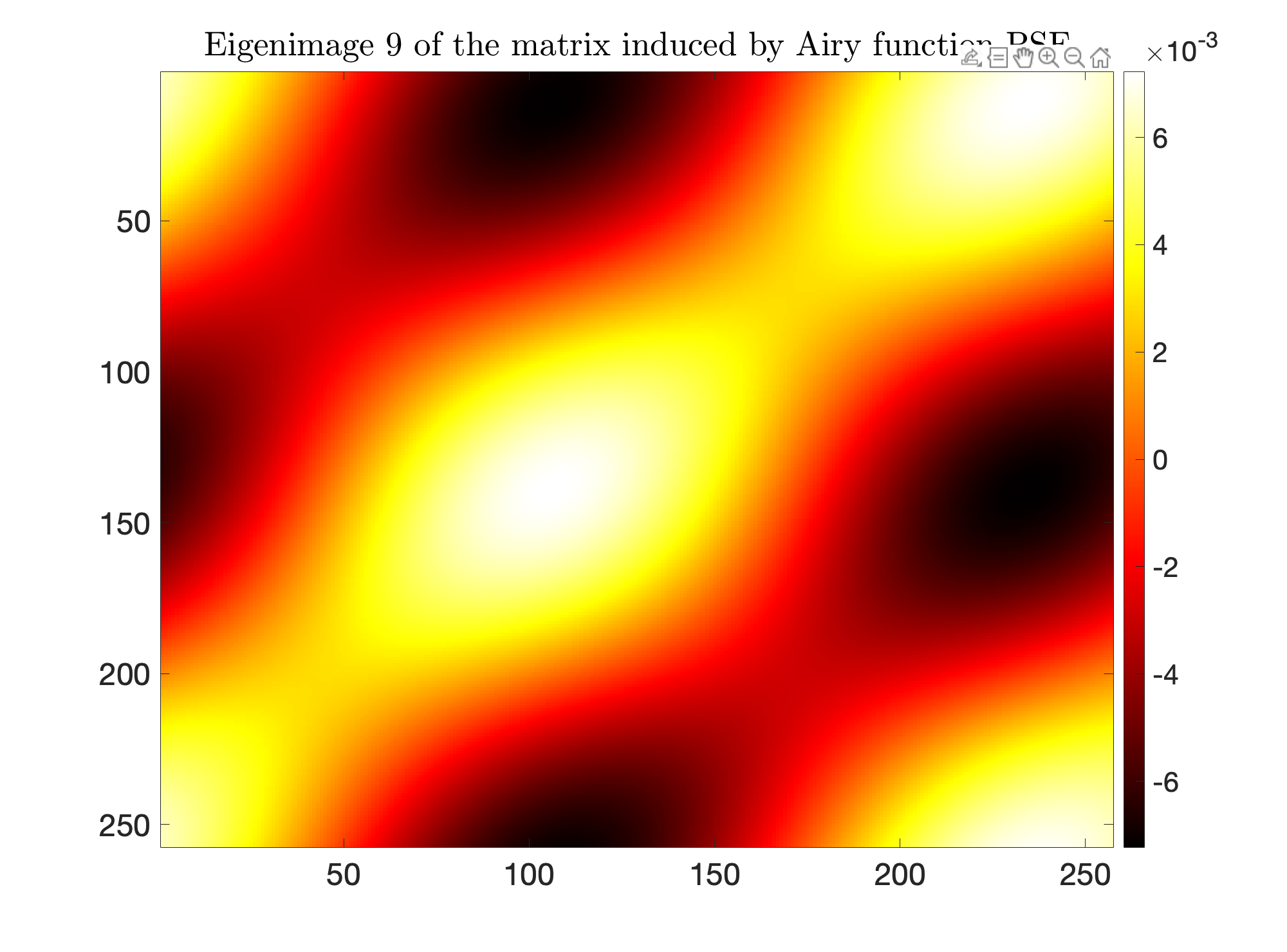}
\end{figure}

It is important to investigate how long it takes to achieve some manner of semiconvergence.  Thus, for each experiment, we track
at which iteration semiconvergence is achieved and also at which iteration the residual-based discrepancy principle is achieved.
For the operator induced by the telescope PSF, we also monitor jumps in eigenvalue magnitude, and these are plotted
separately. This is shown in \Cref{fig.conv-tracker}.  We also show for completeness the first nine eigenvectors (i.e., eigenimages) of the
operator induced by the telescope PSF in \Cref{fig.eigenvectors-airy}.

\begin{figure}
	\centering
	\caption{The image produced by augmented steepest descent at the semiconvergence iteration, raw and with pixel thresholding (all
	pixels of value less than $10$ thresholded to zero\label{fig.xOpt-37-log}), for the experiment augmenting the first 37 eigenvectors.}
	\includegraphics[scale=.15,bb=0 0 1120 840]{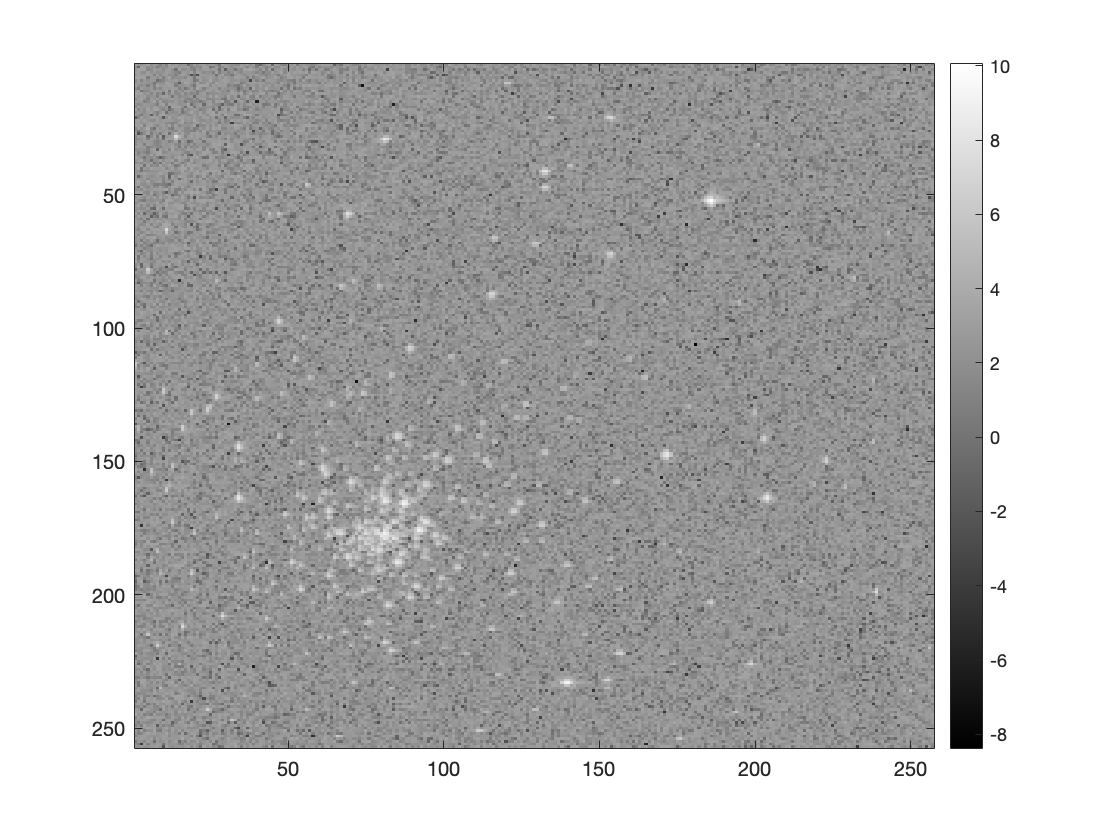}
	\includegraphics[scale=.15,bb=0 0 1120 840]{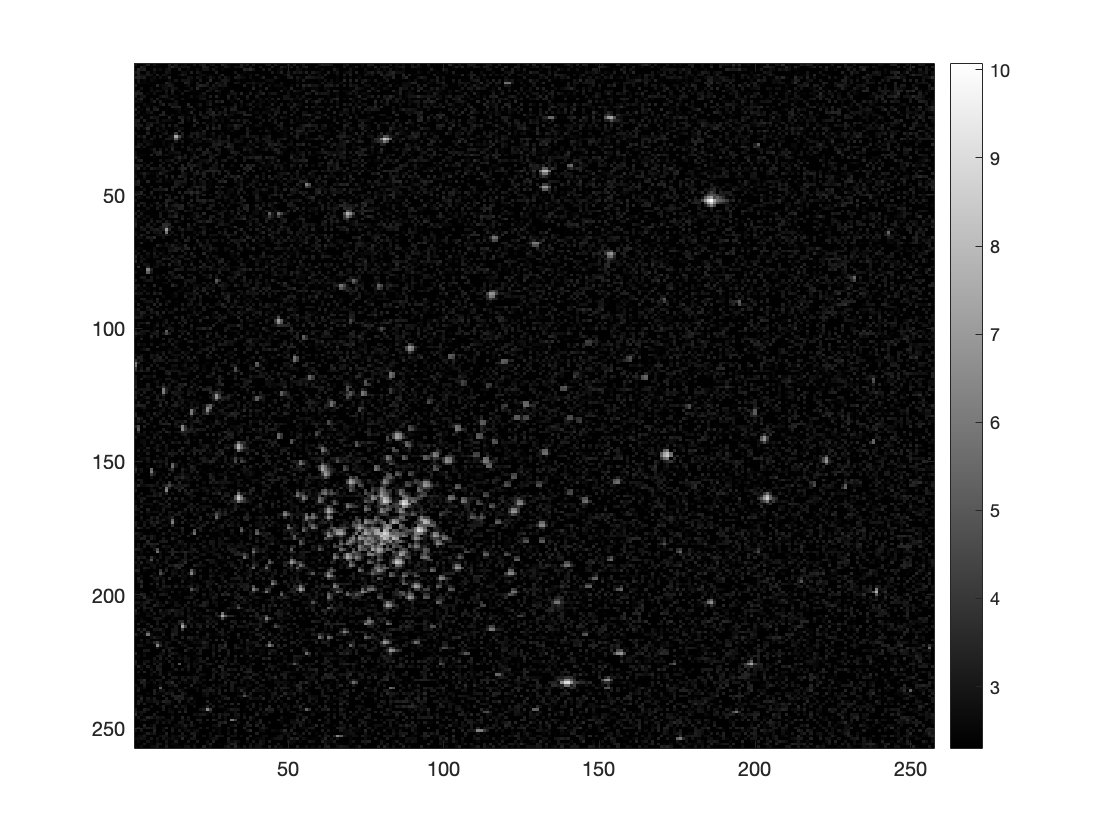}
\end{figure}

Lastly, for the case in which we augmented with the first 37 eigenvectors,
we show the actual image produced at the iteration in which semiconvergence was reached.
We show both the raw image produced by augmented steepest descent, as well as one post-processed by thresholding pixels with values less than $10$
to be zero.  This is shown (again in log plot) in \Cref{fig.xOpt-37-log}.
Augmentation of additional eigenvectors beyond the first 37 produced no additional improvement for this problem.


\section{Conclusions}
In this paper, we have shown that under basic assumptions, the subspace recycling scheme can be combined with any known regularization
scheme, and the resulting augmented scheme is still a regularization.  This opens many possibilities for schemes for developing
new augmented regularization schemes.  We have further demonstrated this by proposing an augmented gradient descent scheme, and
we show how useful recycling easily-computed information within the augmented gradient descent method can be for accelerating
semi-convergence, both in an academic problem and in problems arising from applications in astronomy.

\section*{Acknowledgments}
Victoria Hutterer and Ronny Ramlau are partialy supported by the SFB "Tomography Across the Scales" (Austrian Science Fund Project F68-N36).  The work in this paper was initiated during a visit by the second author to RICAM.
The authors would like to thank Roland Wagner for providing sample PSFs for the last numerical experiment.  Furthermore, the third author thanks
Roland Wagner for his help in using the PSF to induce the adjoint operator needed for applying Landweber-type iterations.
Lastly, the authors thank the referees for their constructive comments. 

\bibliographystyle{siamplain}
\bibliography{references}

\begin{thebibliography}{10}

\bibitem{restore-tools}
{\em Restoretools: An object oriented matlab package for image restoration},
  2012, \url{http://www.mathcs.emory.edu/~nagy/RestoreTools/}.
\newblock Accessed on 28. June 2020.

\bibitem{ASGC-rBiCG-Model-Red.2012}
{\sc K.~Ahuja, E.~de~Sturler, S.~Gugercin, and E.~R. Chang}, {\em Recycling
  {B}i{CG} with an application to model reduction}, SIAM J. Sci. Comput., 34
  (2012), pp.~A1925--A1949, \url{https://doi.org/10.1137/100801500},
  \url{http://dx.doi.org/10.1137/100801500}.

\bibitem{BR-2.2007}
{\sc J.~Baglama and L.~Reichel}, {\em Augmented {GMRES}-type methods},
  Numerical Linear Algebra with Applications, 14 (2007), pp.~337--350,
  \url{https://doi.org/10.1002/nla.518},
  \url{http://dx.doi.org/10.1002/nla.518}.

\bibitem{BR.2007}
{\sc J.~Baglama and L.~Reichel}, {\em Decomposition methods for large linear
  discrete ill-posed problems}, Journal of Computational and Applied
  Mathematics, 198 (2007), pp.~332--343,
  \url{https://doi.org/10.1016/j.cam.2005.09.025},
  \url{http://dx.doi.org/10.1016/j.cam.2005.09.025}.

\bibitem{dSC.2019}
{\sc J.~Chung, E.~de~Sturler, and J.~Jiang}, {\em Hybrid projection methods
  with recycling for inverse problems},  (2020),
  \url{https://arxiv.org/abs/2007.00207}.

\bibitem{deSturler.GCRO.1996}
{\sc E.~{de Sturler}}, {\em Nested {K}rylov methods based on {GCR}}, Journal of
  Computational and Applied Mathematics, 67 (1996), pp.~15--41,
  \url{https://doi.org/10.1016/0377-0427(94)00123-5},
  \url{http://dx.doi.org/10.1016/0377-0427(94)00123-5}.

\bibitem{dSKS.2018}
{\sc E.~de~Sturler, M.~Kilmer, and K.~M. Soodhalter}, {\em Krylov subspace
  augmentation for the solution of shifte systems: a review}.
\newblock in preparation.

\bibitem{demol:1}
{\sc M.~Defrise and C.~D. Mol}, {\em A note on stopping rules for iterative
  methods and filtered svd}, in Inverse Problems: An Interdisciplinary Study,
  P.~C. Sabatier, ed., Academic Press, 1987, pp.~261--268.

\bibitem{DGH.2014}
{\sc Y.~Dong, H.~Garde, and P.~C. Hansen}, {\em R{${}^3$}{GMRES}: including
  prior information in {GMRES}-type methods for discrete inverse problems},
  Electron. Trans. Numer. Anal., 42 (2014), pp.~136--146.

\bibitem{ReRaSoWa}
{\sc L.~Dykes, R.~Ramlau, L.~Reichel, K.~M. Soodhalter, and R.~Wagner}, {\em
  {L}anczos-based fast blind deconvolution methods}, Journal of Computational
  and Applied Mathematics,  (To Appear).

\bibitem{ElVo09}
{\sc B.~Ellerbroek and C.~Vogel}, {\em Inverse problems in astronomical
  adaptive optics}, Inverse Problems, 25 (2009), p.~063001.

\bibitem{RamlauEngl2015}
{\sc H.~Engl and R.~Ramlau}, {\em Regularization of Inverse Problems},
  Springer, New York, 2015, pp.~1233--1241.

\bibitem{Engl_Hanke_Neubauer_1996}
{\sc H.~W. {Engl}, M.~{Hanke}, and A.~{Neubauer}}, {\em {Regularization of
  inverse problems.}}, Dordrecht: Kluwer Academic Publishers, 1996.

\bibitem{Erlangga2008}
{\sc Y.~A. Erlangga and R.~Nabben}, {\em Deflation and balancing
  preconditioners for {Krylov} subspace methods applied to nonsymmetric
  matrices}, SIAM Journal on Matrix Analysis and Applications, 30 (2008),
  pp.~684--699.

\bibitem{Gaul.2014-phd}
{\sc A.~Gaul}, {\em Recycling {K}rylov subspace methods for sequences of linear
  systems: Analysis and applications}, PhD thesis, Technischen Universit\"at
  Berlin, 2014.

\bibitem{GGL.2013}
{\sc A.~Gaul, M.~H. Gutknecht, J.~Liesen, and R.~Nabben}, {\em A framework for
  deflated and augmented {K}rylov subspace methods}, SIAM Journal on Matrix
  Analysis and Applications, 34 (2013), pp.~495--518,
  \url{https://doi.org/10.1137/110820713},
  \url{http://dx.doi.org/10.1137/110820713}.

\bibitem{ir-tools}
{\sc S.~Gazzola, P.~C. Hansen, and J.~G. Nagy}, {\em I{R} {T}ools: a {MATLAB}
  package of iterative regularization methods and large-scale test problems},
  Numer. Algorithms, 81 (2019), pp.~773--811,
  \url{https://doi.org/10.1007/s11075-018-0570-7},
  \url{https://doi-org.libproxy.temple.edu/10.1007/s11075-018-0570-7}.

\bibitem{Gutknecht.AugBiCG.2014}
{\sc M.~H. Gutknecht}, {\em Deflated and augmented {K}rylov subspace methods: a
  framework for deflated {B}i{CG} and related solvers}, SIAM J. Matrix Anal.
  Appl., 35 (2014), pp.~1444--1466, \url{https://doi.org/10.1137/130923087},
  \url{https://doi-org.libproxy.temple.edu/10.1137/130923087}.

\bibitem{Hut18}
{\sc V.~Hutterer and R.~Ramlau}, {\em {Nonlinear wavefront reconstruction
  methods for pyramid sensors using Landweber and Landweber-Kaczmarz
  iteration}}, {Applied Optics}, 57 (2018), pp.~8790--8804.

\bibitem{Hut17}
{\sc V.~Hutterer and R.~Ramlau}, {\em Wavefront reconstruction from
  non-modulated pyramid wavefront sensor data using a singular value type
  expansion}, Inverse Problems, 34 (2018), p.~035002.

\bibitem{HuShaRa19_1}
{\sc V.~Hutterer, R.~Ramlau, and I.~Shatokhina}, {\em Real-time adaptive optics
  with pyramid wavefront sensors: part {I}. a theoretical analysis of the
  pyramid sensor model}, Inverse Problems, 35 (2019), p.~045007,
  \url{https://doi.org/10.1088/1361-6420/ab0656},
  \url{https://doi.org/10.1088%2F1361-6420%2Fab0656}.

\bibitem{ShaHuRam2020}
{\sc R.~R. I.~Shatokhina, V.~Hutterer}, {\em Review on methods for wavefront
  reconstruction from pyramid wavefront sensor data}, Journal of Astronomical
  Telescopes, Instruments and Systems, 6 (2020), p.~010901,
  \url{https://doi.org/10.1117/1.JATIS.6.1.010901}.

\bibitem{Kilmer.deSturler.tomography.2006}
{\sc M.~E. Kilmer and E.~de~Sturler}, {\em Recycling subspace information for
  diffuse optical tomography}, SIAM Journal on Scientific Computing, 27 (2006),
  pp.~2140--2166, \url{https://doi.org/10.1137/040610271},
  \url{http://dx.doi.org/10.1137/040610271}.

\bibitem{Landweber51}
{\sc L.~Landweber}, {\em {An Iteration Formula for Fredholm Integral Equations
  of the First Kind}}, American Journal of Mathematics, 73 (1951),
  pp.~615--624.

\bibitem{LeLouarn_OCTOPUS_04}
{\sc M.~Le~Louarn, C.~V\'{e}rinaud, V.~Korkiakoski, and E.~Fedrigo}, {\em
  {Parallel simulation tools for AO on ELTs}}, in Advancements in Adaptive
  Optics, Proc. SPIE 5490, 2004, pp.~705--712.

\bibitem{LVK06}
{\sc M.~{Le Louarn}, C.~V\'{e}rinaud, V.~Korkiakoski, N.~Hubin, and
  E.~Marchetti}, {\em {Adaptive optics simulations for the European Extremely
  Large Telescope - art. no. 627234}}, in {Advances in Adaptive Optics II, Prs
  1-3}, vol.~6272, {2006}, pp.~{U1048--U1056}.

\bibitem{Louis_1989}
{\sc A.~K. Louis}, {\em {I}nverse und schlecht gestellte {P}robleme}, Teubner
  Studienb{\"u}cher Mathematik, Vieweg+Teubner Verlag, 1989.

\bibitem{morgan.gmresdr}
{\sc R.~B. Morgan}, {\em G{MRES} with deflated restarting}, SIAM J. Sci.
  Comput., 24 (2002), pp.~20--37,
  \url{https://doi.org/10.1137/S1064827599364659},
  \url{https://doi-org.libproxy.temple.edu/10.1137/S1064827599364659}.

\bibitem{Parks.deSturler.GCRODR.2005}
{\sc M.~L. Parks, E.~de~Sturler, G.~Mackey, D.~D. Johnson, and S.~Maiti}, {\em
  Recycling {K}rylov subspaces for sequences of linear systems}, SIAM Journal
  on Scientific Computing, 28 (2006), pp.~1651--1674,
  \url{https://doi.org/10.1137/040607277},
  \url{http://dx.doi.ofrrg/10.1137/040607277}.

\bibitem{parks2016block}
{\sc M.~L. Parks, K.~M. Soodhalter, and D.~B. Szyld}, {\em A block recycled
  gmres method with investigations into aspects of solver performance},
  (2016), \url{https://arxiv.org/abs/1604.01713}.

\bibitem{RaRo12}
{\sc R.~Ramlau and M.~Rosensteiner}, {\em {An efficient solution to the
  atmospheric turbulence tomography problem using Kaczmarz iteration}}, Inverse
  Problems, 28 (2012), p.~095004.

\bibitem{ramlau2020augmented}
{\sc R.~Ramlau and B.~Stadler}, {\em An augmented wavelet reconstructor for
  atmospheric tomography}, 2020, \url{https://arxiv.org/abs/2011.06842}.

\bibitem{S1981}
{\sc Y.~Saad}, {\em Krylov subspace methods for solving large unsymmetric
  linear systems}, Mathematics of Computation, 37 (1981), pp.~105--126,
  \url{https://doi.org/10.2307/2007504},
  \url{http://dx.doi.org/10.2307/2007504}.

\bibitem{Saad.Iter.Meth.Sparse.2003}
{\sc Y.~Saad}, {\em Iterative Methods for Sparse Linear Systems}, SIAM,
  Philadelphia, {S}econd~ed., 2003.

\bibitem{ScherzerBook09}
{\sc O.~Scherzer, M.~Grasmair, H.~Grossauer, M.~Haltmeier, and F.~Lenzen}, {\em
  Variational Methods in Imaging}, Springer, New York, 2009.

\bibitem{soodhalter2020survey}
{\sc K.~M. Soodhalter, E.~de~Sturler, and M.~Kilmer}, {\em A survey of subspace
  recycling iterative methods}, GAMM Mitteilungen,  (2020).
\newblock Applied and numerical linear algebra topical issue.

\bibitem{trefethen.bau}
{\sc L.~N. Trefethen and D.~Bau, III}, {\em Numerical linear algebra}, Society
  for Industrial and Applied Mathematics (SIAM), Philadelphia, PA, 1997,
  \url{https://doi.org/10.1137/1.9780898719574},
  \url{https://doi-org.libproxy.temple.edu/10.1137/1.9780898719574}.

\bibitem{WSG.2007}
{\sc S.~Wang, E.~de~Sturler, and G.~H. Paulino}, {\em Large-scale topology
  optimization using preconditioned {K}rylov subspace methods with recycling},
  International Journal for Numerical Methods in Engineering, 69 (2007),
  pp.~2441--2468, \url{https://doi.org/10.1002/nme.1798},
  \url{http://dx.doi.org/10.1002/nme.1798}.

\bibitem{roland-email}
{\sc R.~Wanger}.
\newblock Personal Communication, Mar. 2020.

\bibitem{YuHeRa13b}
{\sc M.~Yudytskiy, T.~Helin, and R.~Ramlau}, {\em Finite element-wavelet hybrid
  algorithm for atmospheric tomography}, J. Opt. Soc. Am. A, 31 (2014),
  pp.~550--560, \url{https://doi.org/10.1364/JOSAA.31.000550},
  \url{http://josaa.osa.org/abstract.cfm?URI=josaa-31-3-550}.

\end{thebibliography}
\end{document}